\documentclass{article}
\usepackage{graphicx}
\DeclareGraphicsExtensions{.jpg,.pdf}
\usepackage{amssymb,amsmath}
\usepackage[english]{babel}
\usepackage{color}

\begin{document} 


\title{Analysis of a two phase flow model of biofilm spread}

\author{Ana Carpio (Universidad Complutense de Madrid), \\
Gema Duro (Universidad Aut\'onoma de Madrid)}

\date{February 20, 2023}

\maketitle 

{\bf Abstract.} Free boundaries of biofilms advancing on surfaces
evolve according to conservation laws coupled with systems of partial 
differential equations for velocities, pressures and chemicals affecting 
cell behavior.   Thin film approximations lead to complicated 
quasi-stationary systems coupling stationary transport equations and 
compressible Stokes systems with convection-reaction-diffusion 
equations.
We establish existence, uniqueness and stability of solutions of the 
different submodels involved and then obtain well posedness results 
for the full  system.
Our analysis relies on the construction of weak solutions for the steady 
transport equations under sign assumptions and the reformulation 
of the compressible Stokes problem as an elliptic system with enhanced 
regularity properties on the pressure. We need to consider velocity fields 
whose divergence and normal boundary components satisfy sign 
conditions, instead of vanishing as classical results require. 
Applications include the study of cells, biofilms and tissues, where one 
phase is a liquid solution, whereas the other  one is assorted biomass. 
\\
{\bf Keywords.}
Two phase flow, mixture models, thin film approximations,
stationary transport equations,  compressible Stokes equations,
elliptic systems

\section{Introduction}
\label{sec:intro}

Free boundary problems  track the evolution of space regions limited by 
a moving boundary. Films spreading on a surface provide a relevant 
example, with applications in coating, lubrication and biotechnology 
\cite{oron}. In particular, biological films are often described as two phase 
flow mixtures, formed by a biomass phase and a water phase, enclosed 
by a moving boundary \cite{guy,kapellos}. The flow variables are governed 
by sets of coupled conservation laws for mass, momentum and chemical 
species, while the motion of the film boundary is constrained by a
conservation law. Lubrication approximations usually lead to explicit expressions 
for the relevant velocity, pressure and chemical fields, which are used to 
derive a high order nonlinear partial differential equation for the boundary 
dynamics \cite{oron} admitting often explicit self-similar solutions. While this 
approach  yields useful practical information in some regimes, analytical 
studies of the full model are scarce.

\begin{figure}
\centering
\includegraphics[trim=15mm 5cm 15mm 1cm,clip,width=6cm]{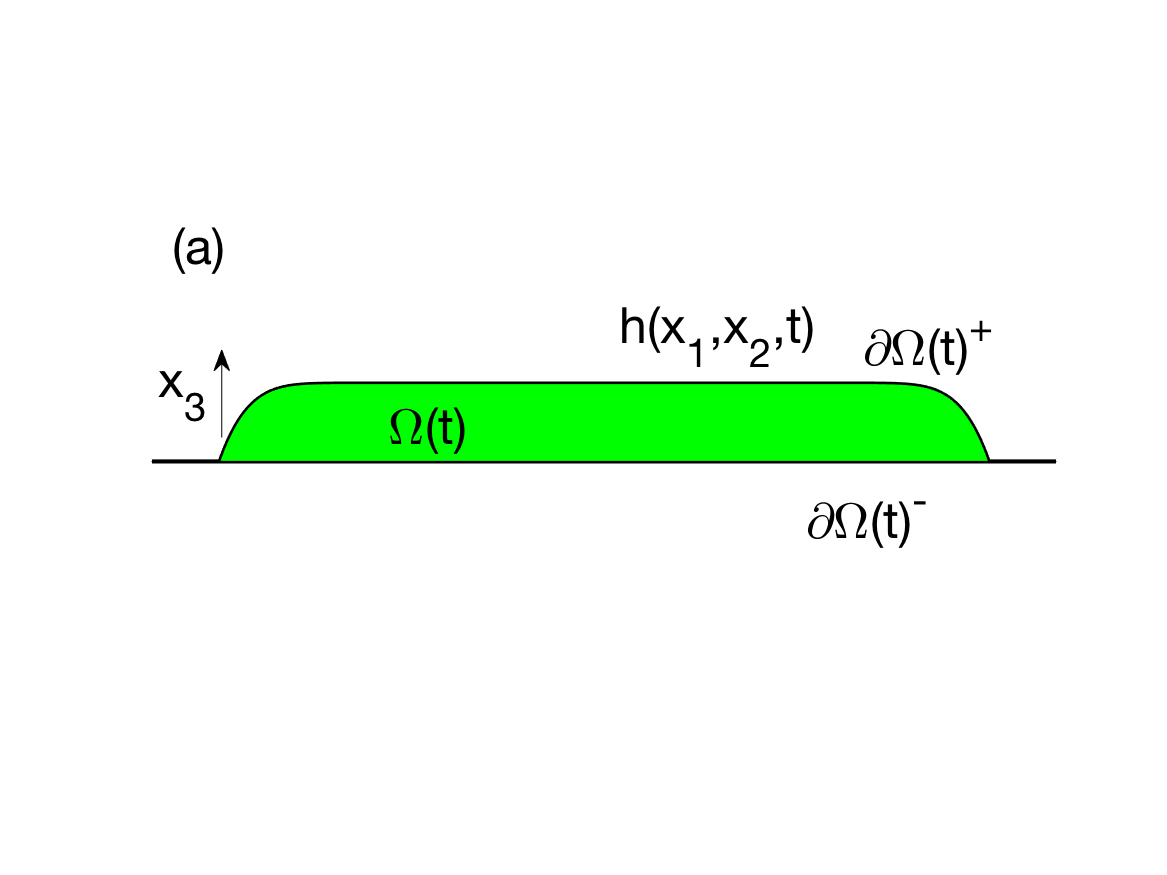}
\includegraphics[trim=15mm 5cm 15mm 1cm,clip,width=6cm]{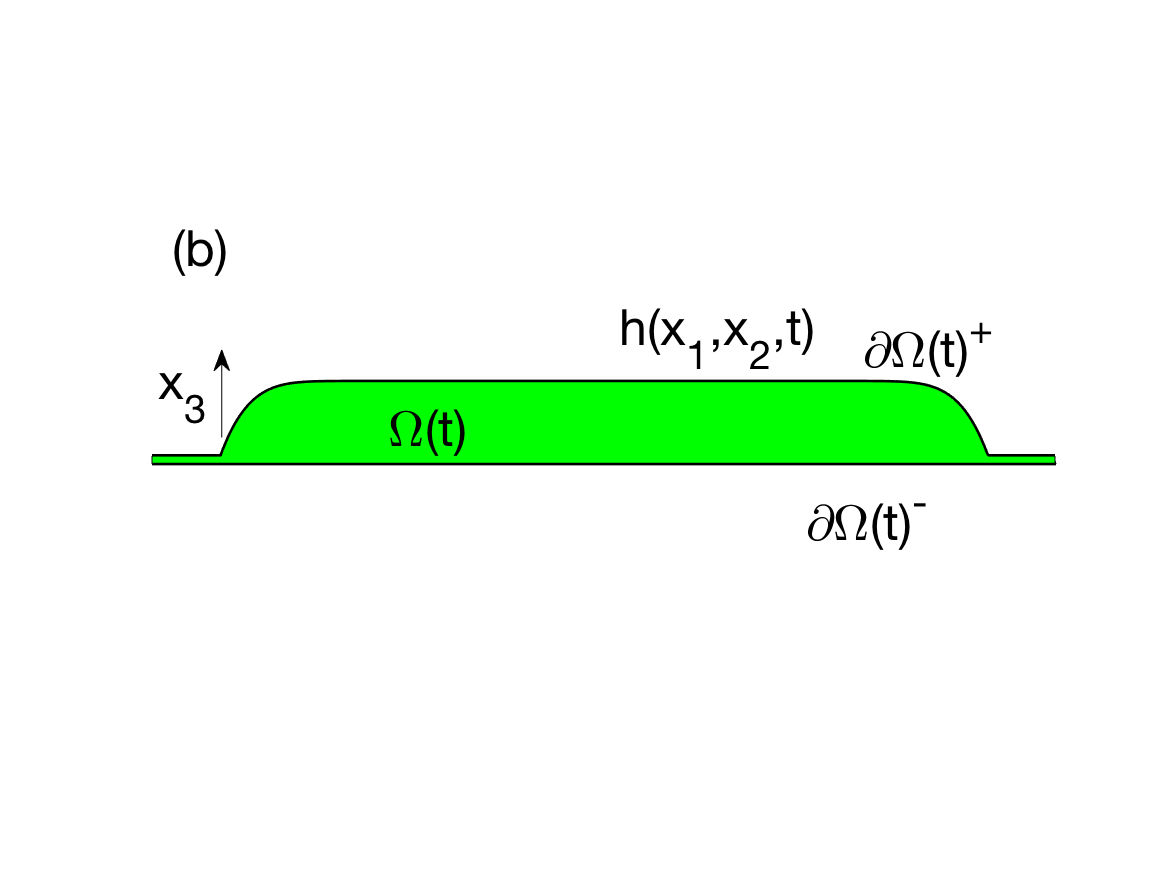}
\caption{Schematic visualization of a biofilm slice.
(a) Biofilm occupying a finite region. The upper boundary joins
the lower boundary at triple points forming angles.
(b) Biofilm with precursor layers around the bulk of the film.}
\label{fig1}
\end{figure}

We aim here to analyze general quasi-stationary systems governing 
mass, velocity, pressure and chemical fields in two phase flow mixture
models. We have in mind applications to the study of cells, biofilms 
and tissues, where one phase is a liquid solution, whereas the other 
one is assorted biomass. To focus the analysis, we consider a specific 
model for bacterial biofilm spread \cite{guy,seminara}:
\begin{eqnarray}
 {\rm div}(\mathbf v_l  \phi_l) =   k_b  {c \over c +K_b} (\phi_l-1),
 \quad  \mathbf x \in \Omega(t), \label{vfractionl_intro} \\ 
 \mu_b \Delta \mathbf v_b +  {\mu_b \over 3}  \nabla {\rm div} 
(\mathbf v_b)  - \nabla (\pi(\phi_b)+ p) =0,  \quad \mathbf x \in \Omega(t),
 \label{velocityb_intro} \\
{\rm div}(\mathbf v_b) = {\rm div}\left(\xi(\phi_l) \nabla p\right),  
\quad  \mathbf x \in \Omega(t), \label{pressure_intro} \\
\mathbf v_l = \mathbf v_b - \eta(\phi_l)  \nabla p,  \;
\phi_l+ \phi_b =1, 
\quad  \mathbf x \in \Omega(t),  \label{velocityl_intro} \\
 - d \Delta c + {\rm div}  (\mathbf v_l c)  
   =  -  k_c \phi_b { c \over c + K_c}, \quad
\mathbf x \in \Omega(t), \label{nutrient_intro}   \\
 {\rm div}(\mathbf v_l  \phi_l+\mathbf v_b  \phi_b) = 0, 
 \quad  \mathbf x \in \Omega(t), \label{boundary_intro} 
\end{eqnarray}
where $\Omega(t)\subset \mathbb R^n$, $n=2,3,$ is the region occupied 
by the biofilm at time $t>0$ (see Figure \ref{fig1}). Here, $\phi_l(\mathbf x,t)$ 
represents the volume fraction of liquid solution and $\phi_b(\mathbf x,t)$ the 
volume fraction of biomass, which move with velocities $\mathbf v_l(\mathbf x,t)$  
and by $\mathbf v_b(\mathbf x,t)$, respectively, under a pressure field $p$.
The parameters $k_b, K_b, k_c, K_c $ and $d$ are positive constants, while
$\pi(\phi_b)$, $\xi(\phi_l)$ and $\eta(\phi_l)$  are known positive functions.
These equations are supplemented with adequate boundary conditions
on $\partial \Omega(t)$.
Notice the total amount of biomass is not constant. It grows thanks to consumption 
of nutrients $c(\mathbf x,t)$. Equations (\ref{vfractionl_intro})-(\ref{nutrient_intro})
define the status of the film variables at each time $t>0$, whereas
the constraint (\ref{boundary_intro}) governs the time dynamics
of the boundary $\partial \Omega(t)$. Section \ref{sec:model} discusses
the model in more detail.

In this paper, we will focus on the solution of system 
(\ref{vfractionl_intro})-(\ref{nutrient_intro}) at a fixed time $t$. 
Equation (\ref{vfractionl_intro}) is a stationary transport equation.
Most existence and regularity results for such equations assume 
$\mathbf v_l \cdot \mathbf n=0$ at the boundary, $\mathbf n$
being the outer unit normal, as well as ${\rm div}(\mathbf
v_l) \in L^\infty$, see \cite{beirao,giraud,novotny}. We will see that 
sign assumptions  allow us to construct weak solutions even when
those two conditions are not satisfied.  
In particular, we assume ${\rm div}(\mathbf v_l)\leq 0$ inside and 
$\mathbf v_l  \cdot \mathbf n \leq 0$ on the boundary, conditions
usually fulfilled by asymptotic and numerical solutions. 
System (\ref{velocityb_intro})-(\ref{pressure_intro}) reminds of
compressible Stokes equations \cite{galdi}. However,
the equation for the pressure leads to more regular pressure fields
here. In fact, (\ref{velocityb_intro})-(\ref{pressure_intro})
can be recast as an elliptic system for the velocity $\mathbf v_b$ 
and pressure $p$ variables with dual data, that is, a right hand
side in Sobolev spaces $W^{-1,q}$. $L^q$ elliptic regularity
for such systems \cite{adn2,kozlov} will be the key  to 
deal with the whole system of equations 
(\ref{vfractionl_intro})-(\ref{pressure_intro}) by an
iterative scheme. The convection-reaction-diffusion equation
(\ref{nutrient_intro}) can be coupled to the scheme provided 
$c$ is uniformly bounded from below by a positive constant.

The paper is organized as follows.
Section \ref{sec:model}  sets up the geometry and presents the
model. We obtain the quasi-static  equations for the magnitudes 
relevant for the  evolution of the film and briefly discuss the 
equations for the  motion of the film boundary.
The rest of the paper is devoted to the analysis of the 
quasi-stationary system.
Sections \ref{sec:transport}, \ref{sec:stokes} and \ref{sec:crd} 
analyze the underlying stationary transport, 
Stokes and convection-reaction-diffusion problems separately.
Section \ref{sec:coupled} proves existence of a solution
of the whole system satisfying a number of stability estimates and 
regularity properties in a fixed domain  by means of an iterative  
scheme. Finally, section \ref{sec:perspectives} summarizes our 
conclusions and discusses the perspectives our work opens to be 
able to handle the equations for the dynamics of the free film
boundaries in general situations.

\section{The two phase flow mixture model}
\label{sec:model}

Assuming that each point in space is occupied by a mixture of
two phases, the composition of the mixture is characterized \cite{guy}
by the volume fractions of both of them. 
We denote by $\phi_l(\mathbf x,t)$ the volume fraction of liquid solution 
and by $\phi_b(\mathbf x,t)$ the volume fraction of biomass (cells
and polymeric matrix), which move with velocities $\mathbf v_l(\mathbf x,t)$ 
and  by $\mathbf v_b(\mathbf x,t)$, respectively. 

\subsection{Conservation laws for the main variables}
\label{sec:conservation}

Conservation of mass for each phase yields the equations
\begin{eqnarray*}
(\rho_l \phi_l)_t + {\rm div}(\mathbf v_l \, \rho_l \phi_l) = J_l,
\\
(\rho_b \phi_b)_t + {\rm div}(\mathbf v_b \, \rho_b \phi_b) =  J_b,
\end{eqnarray*}
where $\rho_l, \rho_b $ denote their densities and $J_l, J_b$
the rate of creation and destruction of each phase. Assuming
there are no external sources, that is, changes only come
from one phase being replaced by the other, we have 
$J_l+J_b=0.$ In the applications considered here, we may set 
$\rho_l =  \rho_b = \rho$ constant. Moreover, creation of biomass 
is the result of cell division due to nutrient consumption. In this
case, the conservation equations become
\begin{eqnarray}
 \phi_{l,t} + {\rm div}(\mathbf v_l  \phi_l) = - r_b(\phi_b,c),
 \label{conservationl} \\
 \phi_{b,t} + {\rm div}(\mathbf v_b  \phi_b) =  r_b(\phi_b,c),
 \label{conservationb}
\end{eqnarray}
where $r_b(\phi_b,c) = k_b \phi_b {c \over c +K_b}$
is the rate of biomass creation,  $k_b, K_b$ being positive
constants. The concentration of nutrients evolves according to
\begin{eqnarray}
{\partial  c \over \partial t} + {\rm div}  (\mathbf v_l c)  
- {\rm div} (d \nabla c)   =  -r_n(\phi_s,c),
\label{nutrient} 
\end{eqnarray}
where $r_n(\phi_b, c) \!=\! k_c \phi_b { c \over c + K_c}$
is the rate of nutrient consumptions. The intake rate $k_c$, the
half saturation $K_c$ and the  diffusivity $d$ are positive 
constants.

Since the volume fractions satisfy 
\begin{eqnarray}
\phi_l(\mathbf x,t)+\phi_b(\mathbf x,t)=1 \label{fraction}
\end{eqnarray}
everywhere, adding equations (\ref{conservationl})-(\ref{conservationb}) 
we conclude that the averaged velocity 
$\mathbf v = \phi_l \mathbf v_l + \phi_b \mathbf v_b$ 
is incompressible
\begin{eqnarray}
{\rm div}(\phi_l \mathbf v_l + \phi_b \mathbf v_b) =
{\rm div}(\mathbf v)=0. 
\label{incompressible}
\end{eqnarray}
However, the fluid velocities in each phase are not divergence
free in principle, that is, the fluid phases may be compressible.
Notice that (\ref{incompressible}) implies that $\phi_l + \phi_b$ is
constant in time. If $\phi_l(0) + \phi_b(0)=1$, then 
$\phi_l(t) + \phi_b(t)=1$ for $t>0$.

In most cellular samples the velocities $\mathbf v_l$ and $\mathbf v_b$ 
are small enough for inertial forces to be negligible. Then, the velocities 
are determined from a balance of  forces in the two fluids 
\cite{guy,kapellos,lanir}
\begin{eqnarray}
{\rm div} \mathbf T_l  + F_{b\rightarrow l} + \phi_l \nabla \pi_l = 0,
\label{momentuml} \\
{\rm div} \mathbf T_b +F_{l\rightarrow b}  + \phi_b \nabla \pi_b = 0,
\label{momentumb}
\end{eqnarray}
where the forces associated to chemical potentials satisfy
$\phi_b\nabla \pi_b  + \phi_l \nabla \pi_l = 0$, and the interaction
forces of one phase on the other through the interfaces too, 
$F_{l\rightarrow b}  + F_{b\rightarrow l}=0$.
The stress tensors for the viscous fluid and the aqueous liquid phases 
are given by
\begin{eqnarray*}
\mathbf T_l = - \phi_l p_l \, \mathbf I,  \label{stressl} \\
\mathbf T_b = - \phi_b p_b \, \mathbf I + \boldsymbol \sigma_b, \quad
\boldsymbol \sigma_b = \mu_b (\nabla \mathbf v_b + \nabla \mathbf v_b^t)
- 2 {\mu_b\over 3} {\rm div}(\mathbf v_b) \, \mathbf I, \label{stressb}
\end{eqnarray*}
where $\mu_b$ is the shear viscosity and $\mathbf I$ the identity.
The interaction forces are described by the constitutive law \cite{guy,kapellos}
\begin{eqnarray*}
F_{l\rightarrow b} =  K(\phi_l) (\mathbf v_l - \mathbf v_b) + 
p_{lb} \nabla \phi_b = - F_{b\rightarrow l}.
\end{eqnarray*}
To simplify, we take the intraphase pressures $p_l$, $p_b$ and
the interphase pressure $p_{lb}$ equal, that is, $p_l=p_b=p_{lb}=P$.
Otherwise, additional constitutive laws would be needed for them
\cite{guy}.
Using these expressions, equations (\ref{momentuml})-(\ref{momentumb})
become
\begin{eqnarray}
- \phi_l \nabla P + K(\phi_l) (\mathbf v_b - \mathbf v_l)  
+ \phi_l  \nabla \pi_l = 0.
\label{motionl1} \\
{\rm div} \,  \boldsymbol \sigma_b   - \phi_b \nabla P 
+ K(\phi_l) (\mathbf v_l - \mathbf v_b) +   \phi_b  \nabla \pi_b  = 0,
\label{motionb1}.
\end{eqnarray}
Setting now $p= P-\pi_l$  and adding (\ref{motionl1})-(\ref{motionb1})
we find the equations
\begin{eqnarray}
 K(\phi_l) (\mathbf v_b - \mathbf v_l) - \phi_l  \nabla p = 0,
\label{motionl} \\
\mu_b \Delta \mathbf v_b +  {\mu_b \over 3}  \nabla {\rm div} 
(\mathbf v_b)  - \nabla (\pi_l + p)  = 0.
\label{motionb}
\end{eqnarray}
The  pressure variable $p$ allows us to satisfy the incompressibility
condition (\ref{incompressible}). Combining (\ref{incompressible})
with (\ref{fraction}) and (\ref{motionl}), we get an equation for
the pressure
\begin{eqnarray}
{\rm div}(\mathbf v_b) = {\rm div}\left({\phi_l^2 \over K(\phi_l)} \nabla p
\right).
\label{pressure}
\end{eqnarray}
Typically, $K(\phi_l) = {\mu_l \over \zeta(\phi_l)}$ with $\zeta(\phi_l)
\sim \zeta_{\infty}$  and $\pi_l = \pi(\phi_b) = \Pi \phi_b$, $\Pi, \mu_l>0$
\cite{seminara}.

\subsection{Quasi-stationary approximation and boundary conditions}
\label{sec:final}

The final set of equations is given by (\ref{conservationl}), (\ref{nutrient}), 
(\ref{fraction}), (\ref{incompressible}), (\ref{motionl}), (\ref{motionb}), 
(\ref{pressure}),  posed in the region  $\Omega(t)$ occupied by the biological 
sample. Figure \ref{fig1} represents a two dimensional slice of a three 
dimensional biofilm. We can work with three dimensional objects or
two dimensional slices to simplify. In any case, we have an horizontal 
bottom boundary $\partial \Omega(t)^-$ and an upper moving boundary
$\partial \Omega(t)^+$, and two standard configurations. In one of them
the borders of $\partial \Omega(t)^+$ join $\partial \Omega(t)^-$ at triple
points, see Fig. \ref{fig1}(a). In the other $\partial \Omega(t)^+$ decays
at the biofilm edges to form very thin precursor layers \cite{degennes},
see Fig. \ref{fig1}(b).

In the biological  applications we have in mind we may neglect the time 
derivatives in (\ref{conservationl}) and (\ref{nutrient}).   
The resulting approximated system in dimensionless form is given by 
(\ref{vfractionl_intro}), (\ref{velocityb_intro}), (\ref{pressure_intro}), 
(\ref{velocityl_intro}), (\ref{nutrient_intro})
where $k_b, K_b, k_c, K_c, d$ are positive constants  \footnote{
For simplicity, we keep the same symbols as in the previous dimensional 
equations. However, they now represent dimensionless combinations
of the dimensional parameters and the chosen scales}. We set
$\pi(\phi_b) = \Pi \phi_b$, $\xi(\phi_l) \sim {\mu_l \over \zeta_\infty} \phi_l^2
\sim {\mu_l \over \zeta_\infty} \phi_\infty^2$
and $\eta(\phi_l) \sim  {\xi(\phi_l) \over \phi_l}  \sim
{\mu_l \over \zeta_\infty} \phi_\infty$, with $\Pi, \mu_l, \zeta_\infty >0,$ 
$\phi_\infty \in (0,1)$.
With these choices, taking the divergence of  equation (\ref{velocityb_intro}), 
we find the  additional relation
\begin{eqnarray}
{4 \mu_b \over 3}  \Delta  {\rm div} (\mathbf v_b) = \Pi \Delta \phi_b 
 + \Delta p,  \quad \mathbf x \in \Omega(t).
\label{divergence_fin} 
\end{eqnarray}


As for the boundary conditions, 
\begin{eqnarray}
\left( \mu_b (\nabla \mathbf v_b + \nabla \mathbf v_b^t) - (
2 {\mu_b \over 3} {\rm div}(\mathbf v_b) + p + \Pi \phi_b) \mathbf I \right)
\cdot \mathbf n = \mathbf t_{\rm ext}, & \;
\mathbf x \in \partial \Omega(t)^+, \label{bc_velocity+} \\
\mathbf v_b = 0, & \; \mathbf x \in \partial \Omega(t)^-, 
\label{bc_velocity-} \\
p = p_{\rm ext} - \pi_{\rm ext}, & \mathbf x \in \partial \Omega(t),
\label{bc_pressure+-} \\
c = c_0,  &\, \mathbf x \in \partial \Omega(t)^-,  \label{bc_concentration-} \\
{\partial c \over \partial \mathbf n} = 0,  & \, 
\mathbf x \in \partial \Omega(t)^+, \label{bc_concentration+}
\end{eqnarray}
where $\mathbf t_{\rm ext}$, $p_{\rm ext}$, $\pi_{\rm ext}$ and
$c_0$ are given external tensions, pressures and concentrations.

Since all unknowns depend on time through the motion of the interface, 
the final problem is quasi-stationary.

\subsection{Dynamics of the moving interface}
\label{sec:interface}

Assuming the boundary $\partial \Omega(t)$ is defined by a
surface $x_3=h(x_1,x_2,t)$, we can obtain an equation for its
dynamics integrating equation (\ref{incompressible}) in the $x_3$
direction to find
\begin{eqnarray*}
\begin{array}{l}
\int_{0}^h {\partial (\mathbf v \cdot \hat{\mathbf x}_1) \over \partial x_1} \, dx_3 
+ \int_{0}^h {\partial (\mathbf v \cdot \hat{\mathbf x}_2) \over \partial x_2} \, dx_3 
+ \int_{0}^h {\partial (\mathbf v \cdot \hat{\mathbf x}_3) \over \partial x_3} \, dx_3 
=0,
\end{array}
\end{eqnarray*}
where $\hat{\mathbf x}_1$, $\hat{\mathbf x}_2$ and $\hat{\mathbf x}_3$ are the  
unit vectors in the cartesian coordinate directions. By Leibniz's rule:
\begin{eqnarray*}
\begin{array}{l}
\int_{0}^h {\partial (\mathbf v \cdot \hat{\mathbf x}_i) \over \partial x_i} \, dx_3
= {\partial \over \partial x_i} \left[
\int_{0}^h (\mathbf v \cdot \hat{\mathbf x}_i) \, dx_3 \right]
- \mathbf v \cdot \hat{\mathbf x}_i\big|_h {\partial h \over \partial x_i}, 
\quad i=1,2.
\end{array}
\end{eqnarray*}
Thus
\begin{eqnarray}
\begin{array}{l}
{\partial \over \partial x_1} \left[
\int_{0}^h (\mathbf v \cdot \hat{\mathbf x}_1) \, dx_3 \right]
+ 
{\partial \over \partial x_2} \left[
\int_{0}^h (\mathbf v \cdot \hat{\mathbf x}_2) \, dx_3 \right]
\\  
- \mathbf v \cdot \hat{\mathbf x}_1\big|_h {\partial h \over \partial x_1}
- \mathbf v \cdot \hat{\mathbf x}_2\big|_h {\partial h \over \partial x_2}
+ \mathbf v \cdot \hat{\mathbf x}_3\big|_h  =
\mathbf v \cdot \hat{\mathbf x}_3\big|_{x_3=0}.
\end{array} 
\label{height1}
\end{eqnarray}
Next, we differentiate $x_3(t)=h(x_1(t),x_2(t),t)$ with respect to time and
use $\mathbf v \cdot \hat{\mathbf x}_i = {dx_i \over dt}$, $i=1,2,3,$  to get
\begin{eqnarray*}
\begin{array}{lll}
\mathbf v \cdot x_3 \big|_h = {d x_3 \over dt} = {d \over dt} h(x_1(t),x_2(t),t) 
&=& {\partial h \over \partial t} + {\partial h \over \partial x_1} {d x_1 \over dt}
+ {\partial h \over \partial x_2} {d x_2 \over dt}  \\[1ex]
&=& {\partial h \over \partial t} + 
\mathbf v \cdot x_1 \big|_h  {\partial h \over \partial x_1} 
+ \mathbf v \cdot x_2 \big|_h {\partial h \over \partial x_2}.
\end{array}
\end{eqnarray*}
Inserting this identity in (\ref{height1}) we find the equation
\begin{eqnarray}
{\partial h \over \partial t } +
{\partial \over \partial x_1} \left[
\int_{0}^h (\mathbf v \cdot \hat{\mathbf x}_1) \, dx_3 \right]
+ {\partial \over \partial x_2} \left[
\int_{0}^h (\mathbf v \cdot \hat{\mathbf x}_2) \, dx_3 \right]  =
\mathbf v \cdot \hat{\mathbf x}_3\big|_{x_3=0},
\label{height2}
\end{eqnarray}
where $\mathbf v \cdot \hat{\mathbf x}_i = v_i = v_{b,i} 
 - \eta(\phi_l)  {\partial p \over \partial x_i}, i=1,2,3,$
or, equivalently,
\begin{eqnarray}
{\partial h \over \partial t } 
+ v_1 \big|_{x_3=h} {\partial h \over \partial x_1} 
+ v_2 \big|_{x_3=h} {\partial h \over \partial x_2}
+ \int_{0}^h {\partial v_1  \over \partial x_1} \, dx_3 
+ \int_{0}^h {\partial v_2  \over \partial x_2} \, dx_3  
= v_3\big|_{x_3=0}.
\label{height3}
\end{eqnarray}
These equations hold for $(x_1,x_2) \in \partial \Omega(t)^-$,
that is, for points belonging to the bottom boundary. At the
edges we must apply boundary conditions. To that purpose,
the configuration  represented in Fig. \ref{fig1}(b) is easier
to  handle, we can enforce zero Neumann boundary conditions
or an asymptotically constant value.
The equations can be applied to a two dimensional slice
by just dropping the $x_2$ variable from the equations and
working with $x_1$ and $x_3$. \\

The two phase flow problem under study consists of equation (\ref{height2}) 
for the motion of the interface $h$ defining $\Omega(t)$, coupled to the set of 
quasi-stationary equations for the velocities $\mathbf v_b$,$\mathbf v_l$, the 
pressure $p$, the volume fractions $\phi_b$,$\phi_l$ and the nutrient 
concentration $c$ set in $\Omega(t)$. 
Next, we consider each of these equations separately and explain how to 
construct solutions for them in a fixed domain $\Omega$.

\section{The stationary transport problem}
\label{sec:transport}

For sign reasons, we choose to work with the scalar equation for 
$\phi_l = 1 - \phi_b$
\begin{eqnarray}
 {\rm div}(\mathbf v_l  \phi_l) =   - k_b \phi_b {c \over c + K_b} = 
- k_b {c \over c + K_b} + k_b {c \over c + K_b} \phi_l,
 \quad   \mathbf x \in \Omega,   \label{vfractionl_st}
\end{eqnarray}
which is equivalent to a similar equation for $\phi_b$ thanks to 
(\ref{incompressible}). Here, $k_b$ and $K_b$ are positive constants 
and $c, \mathbf v_l$ are known functions.  

In practice, to study the full coupled model one often sets
${c \over c+K_b} = g_\infty >0$  in (\ref{vfractionl_st}),
leading to the scalar equation:
\begin{eqnarray}
 {\rm div}(\mathbf v_l  \phi_l) =   - k_b g_\infty \phi_b  = 
- k_b g_\infty + k_b g_\infty \phi_l,
 \quad   \mathbf x \in \Omega.   \label{vfractionl_stinfty}
\end{eqnarray}

Unlike stationary transport problems often studied \cite{beirao,giraud},
these equations contain the term  ${\rm div}(\mathbf v_l)  \phi_l$ and
we cannot assume $\mathbf v_l \cdot \mathbf n =0$. From numerical
simulations and asymptotic calculations \cite{entropy,seminara},
we expect ${\rm div}(\mathbf v_l) \leq 0$ and $\mathbf v_l \cdot \mathbf n 
\leq 0$. We establish next an existence and regularity theory for stationary
transport problems in divergence form 
\begin{eqnarray}
 {\rm div}(-\mathbf v_l(\mathbf x)  \phi_l) + a(\mathbf x) \phi_l = 
 a(\mathbf x), \quad   \mathbf x \in \Omega,   \label{vfractionl_gen}
\end{eqnarray}
under these conditions \footnote{The advection operator in equation 
 ${\rm div}(- \mathbf v_l  \phi) + a \phi  = g$ 
only needs boundary conditions for $\phi$ on 
$\Gamma = \{ \mathbf x \in \partial \Omega \, | \, -\mathbf v_l  
\cdot \mathbf n < 0 \} \subset \partial \Omega,$ 
see \cite{friedrichs}.}. \\


{\bf Theorem 3.1 (Existence).} {\it Let $\Omega \subset \mathbb R^n$, $n=2,3$ 
be an  open, bounded set, with $C^1$ boundary $\partial \Omega$. Let 
$\mathbf v_l \in [H^1(\Omega) \cap C(\overline{\Omega})]^n$ be such that  
${\rm div}(\mathbf v_l)\leq 0$ in  $\Omega$ and $\mathbf v_l \cdot \mathbf 
n \leq 0$  on $\partial \Omega$. Let $a \in L^\infty(\Omega)$  be a strictly 
positive function  bounded from below by a positive constant $a_{\rm min}$, 
that is, $a \geq a_{\rm min}>0$ in $\Omega$.
Then, there exists a weak solution $\phi_l \in L^2(\Omega)$ of 
(\ref{vfractionl_gen}) in the sense of distributions. Moreover,
$ 0 \leq \phi \leq 1$ a.e. on $\Omega$ and $\phi$ cannot vanish in a
set of positive measure.}

{\bf Proof.} {\it Existence.} 
Following \cite{beirao}, for each $\varepsilon >0$, let 
$\phi_\varepsilon \in  H^1(\Omega)$ be  the solution of
\begin{eqnarray}
- \varepsilon \Delta  \phi _\varepsilon  - {\rm div}(\mathbf v_l  \phi_\varepsilon)
+ a \phi_\varepsilon = a \; {\rm in \,} \Omega, \quad 
{\partial  \phi_\varepsilon  \over \partial \mathbf n } = 0  \; {\rm on \,} 
\partial \Omega.
\label{eq_phie}
\end{eqnarray}
In variational form, the equation reads
\begin{eqnarray*} b(\phi_\epsilon, w) =
\varepsilon \int_\Omega \nabla \phi _\varepsilon \cdot \nabla w \, d \mathbf x
+ \int_\Omega \mathbf v_l  \cdot \phi _\varepsilon \nabla w \, d \mathbf x
- \int_{\partial \Omega} \phi_\varepsilon  w \, \mathbf v_l \cdot \mathbf n
d S_{\mathbf x} \\
+ \int_\Omega a \, \phi_\varepsilon w \, d \mathbf x  
= \int_\Omega a \, w \, d\mathbf x  = L(w),
\end{eqnarray*}
for $w \in H^1(\Omega)$,  defined on $\partial \Omega$ as 
$L^2(\partial \Omega)$ functions in the sense of traces \cite{brezis}. 
The bilinear form $b(\varphi, w)$ is continuous on $H^1(\Omega)$ 
\cite{raviart2}, while the linear form $L$ is continuous on $L^2(\Omega).$ 

Since ${\rm div}(\mathbf v_l) \leq 0$, $\mathbf v_l \cdot  \mathbf n \leq 0$ 
and $a>a_{\rm min}$, the bilinear form $b$ is also coercive in 
$H^1(\Omega)$. Notice that
\begin{eqnarray*}
\int_\Omega \hskip -1mm  \mathbf v_l \cdot \phi_\varepsilon  \nabla 
\phi_\varepsilon  d \mathbf x
= {1\over 2}\int_\Omega \hskip -1mm  \mathbf v_l  \cdot \nabla 
|\phi_\varepsilon|^2  d \mathbf x =  
{1\over 2}\int_{\partial \Omega} \hskip -1mm |\phi_\varepsilon|^2 
\mathbf v_l \cdot \mathbf n d \mathbf x
 - {1\over 2} \int_\Omega \hskip -1mm {\rm div}(\mathbf v_l) 
|\phi_\varepsilon|^2   d \mathbf x.
\end{eqnarray*}
The positive term $ - \int_\Omega {\rm div}(\mathbf v_l) |\phi_\varepsilon|^2 
d \mathbf x $ is finite because $|\phi_\varepsilon|^2 \in L^2(\Omega)$,
that is,  $\phi_\varepsilon \in L^4(\Omega)$ thanks to Sobolev embeedings
\cite{adams,brezis}. Moreover, the bilinear form
$\varepsilon \int_\Omega \nabla \phi \cdot \nabla w \, d \mathbf x
+  \int_\Omega a \phi w \, d \mathbf x$
is coercive in $H^1(\Omega)$, see \cite{raviart2}.
Thus, by Lax Milgram's theorem, we have a  unique solution 
$\phi_\varepsilon \in H^1(\Omega)$, see \cite{brezis}. 

Setting $w=\phi _\varepsilon$,  Young's inequality \cite{brezis} implies
\begin{eqnarray*} 
0 \leq \varepsilon \int_\Omega |\nabla \phi _\varepsilon|^2 \, d \mathbf x
- {1\over 2} \int_{\partial \Omega}  |\phi_\varepsilon|^2  \mathbf v_l \cdot 
\mathbf n \, d S_{\mathbf x} +  \int_\Omega \left[ - {1\over 2} 
{\rm div}(\mathbf v_l)  + a \right] |\phi_\varepsilon|^2  \, d \mathbf x  \\ 
= \int_\Omega a \phi_\varepsilon  \, d\mathbf x   
\leq \| a \|_{L^2} \left( \int_\Omega  |\phi_\varepsilon|^2 \right)^{1/2}.
\end{eqnarray*}
Thus, $ a_{\rm min}\| \phi_\varepsilon \|_{L^2} \leq \| a \|_{L^2}$. Then,
each of the additional positive terms in the left hand side of the above 
inequality are uniformly bounded in terms of $\| a \|_{L^2}$.
Therefore, we can extract a subsequence $\phi_{\varepsilon'}$ such that 
$ \phi_{\varepsilon'}$ tends weakly in $L^2(\Omega)$ to a limit $\phi$,
and $\varepsilon \nabla \phi_\varepsilon$ tends strongly to zero.
Setting now $w \in C_c^\infty(\Omega)$ in the variational formulation,
and passing to the limit \cite{aml,profiles}
we see that $\phi$ is a solution of (\ref{vfractionl_st}) in the sense
of distributions. 
The variational formulation holds with $\epsilon =0$ and the
boundary integral replaced by the duality product 
$_{H^{-1/2}(\partial \Omega)}<\phi \, \mathbf v_l \cdot 
\mathbf n, w>_{H^{1/2}(\partial \Omega)}$ for $w\in H^1(\Omega)$.

{\it $L^\infty$ bounds. }
The functions $\psi_\varepsilon =  \phi_\varepsilon - M$ satisfy
\begin{eqnarray*} 
\varepsilon \int_\Omega \nabla \psi _\varepsilon \cdot \nabla w \, d \mathbf x
+ \int_\Omega \mathbf v_l \cdot \psi _\varepsilon \nabla w \, d \mathbf x
- \int_{\partial \Omega} \psi_\varepsilon  w \, \mathbf v_l \cdot \mathbf n 
d S_{\mathbf x} + \\
\int_\Omega a  \psi_\varepsilon w \, d \mathbf x  
= \int_\Omega \left[{\rm div}(\mathbf v_l) M + a (1-M) 
\right] \, w \, d\mathbf x.
\end{eqnarray*}
Setting $M= 1$ and  $w = \psi_\varepsilon^+$, we get
\begin{eqnarray*} 
\varepsilon \int_\Omega |\nabla \psi _\varepsilon^+|^2 \, d \mathbf x - 
{1\over 2} \int_{\partial \Omega}  |\psi_\varepsilon^+|^2 \, 
\mathbf v_l \cdot \mathbf nd S_{\mathbf x} 
+  \int_\Omega \left[ - {1\over 2} {\rm div}(\mathbf v_l)  +
a \right] |\psi_\varepsilon^+|^2  \, d \mathbf x    \\ 
= \int_\Omega  {\rm div}(\mathbf v_l) \psi_\varepsilon^+ d\mathbf x \leq 0.
\end{eqnarray*}
Thus, $\psi_\varepsilon^+=0$ and $ \phi_\varepsilon \leq 1.$ Similarly,
we set  $\psi_\varepsilon = M - \phi_\varepsilon$ and $M=0$ to find
\begin{eqnarray*} 
\varepsilon \int_\Omega |\nabla \psi _\varepsilon^+|^2 \, d \mathbf x - 
{1\over 2} \int_{\partial \Omega} (\mathbf v_l \cdot \mathbf n) |\psi_\varepsilon^+|^2 
\, d S_{\mathbf x} +  \int_\Omega \left[ - {1\over 2} {\rm div}(\mathbf v_l)  +
a \right] |\psi_\varepsilon^+|^2  \, d \mathbf x    \\ 
= - \int_\Omega  a  \psi_\varepsilon^+ d\mathbf x \leq 0.
\end{eqnarray*}
Thus, $\psi_\varepsilon^+=0$ and $ \phi_\varepsilon \geq 0.$
Any weak limit $ \phi$ in $L^2$ inherits these two properties in the region where
$  - {1\over 2} {\rm div}(\mathbf v_l)  + a   \neq 0.$ 
Since $a$ is strictly positive, this happens a.e. in the whole set $\Omega$.
Moreover, if $\phi$ vanish in a set of positive measure, then
(\ref{vfractionl_gen}) implies that $a$ vanishes in the same set, which is
impossible.
$\square$ \\


{\bf Theorem 3.2 ($L^2$ Regularity and uniqueness).} {\it Under
the assumptions of Theorem 3.1, if $\partial \Omega$ has $C^2$ 
regularity
\footnote{We could handle piecewise $C^2$ domains and convex 
Lipschitz domains using results in Sobolev spaces and elliptic 
regularity in such domains \cite{grisvard, necas}.}, 
$\nabla a \in L^2(\Omega)$, $\mathbf v_l \in H^2(\Omega)$ 
and  $ \nabla \mathbf v_l \in [L^{\infty}(\Omega)]^{n^2}$ with  
$\|\nabla \mathbf v_l\|_{[L^{\infty}]^{n^2}}$ small enough 
compared to $a_{\rm min}$, then $\|\nabla \phi \|_{[L^2]^n}$ is 
bounded from above by a constant depending on 
$\|\nabla {\rm div}(\mathbf v_l)\|_{[L^2]^n}$, $a_{\rm min}$
and $\|\nabla a \|_{[L^2]^n}$.
The solution $u \in H^1(\Omega)$ is unique.
}

{\bf Proof.} {\it Regularity.}
Elliptic regularity applied to system (\ref{eq_phie}) implies that its solution 
$\phi_\varepsilon \in H^2(\Omega)$ \cite{adn,gilbart}.
Multiplying equation (\ref{eq_phie}) by $\Delta \phi_\varepsilon$ and 
integrating over $\Omega$ we find
\begin{eqnarray*}
-  \varepsilon \int_\Omega |\Delta \phi _\varepsilon|^2  d \mathbf x
- \int_\Omega \!\! \mathbf v_b \!\cdot\! \nabla \phi_\varepsilon \Delta 
\phi _\varepsilon d \mathbf x
+ \int_\Omega \!\! \left[ - {\rm div}(\mathbf v_b) + a
\right] \phi_\varepsilon  \Delta \phi_\varepsilon  d \mathbf x  
= \int_\Omega a \Delta \phi _\varepsilon d \mathbf x.
\end{eqnarray*}
Let us rewrite the second integral term. To simplify, we use the summation
convention, that is, sum over repeated indexes is intended. Integrating by
parts we get
\begin{eqnarray*}\begin{array}{l}
 - \int_\Omega \mathbf v_l \cdot \nabla \phi_\varepsilon  \Delta \phi_\varepsilon  
d \mathbf x =  
- \int_\Omega v_{l,j} \phi_{\varepsilon, x_j}  \phi_{\varepsilon, x_k x_k} d \mathbf x 
 \\[1.5ex]
 = \int_\Omega v_{l,j,x_k} \phi_{\varepsilon, x_j}  \phi_{\varepsilon, x_k} 
d \mathbf x + \int_\Omega v_{l,j} \phi_{\varepsilon, x_j x_k}  \phi_{\varepsilon, x_k} 
d \mathbf x  \\[1.5ex]
 = \int_\Omega v_{l,j,x_k} \phi_{\varepsilon, x_j}  \phi_{\varepsilon, x_k}  d \mathbf x
- {1\over 2} \int_\Omega {\rm div}(\mathbf v_l) |\nabla \phi_{\varepsilon}|^2  d \mathbf x
+ {1\over 2} \int_{\partial \Omega}  |\nabla \phi_{\varepsilon}|^2 \mathbf v_l \cdot 
\mathbf n \, d S_{\mathbf x}.  
\end{array}\end{eqnarray*}
Integrating by parts again and using the boundary condition, the remaining 
terms give
\begin{eqnarray*}
- \int_\Omega {\rm div}(\mathbf v_l) \phi_\varepsilon  \Delta \phi_\varepsilon  
d \mathbf x
=  \int_\Omega {\rm div}(\mathbf v_l)  |\nabla \phi_\varepsilon|^2  d \mathbf x
+ \int_\Omega \nabla {\rm div}(\mathbf v_l) \cdot \phi_\varepsilon  
\nabla \phi_\varepsilon   d \mathbf x, \\
\int_\Omega a \phi_\varepsilon  \Delta \phi_\varepsilon  
d \mathbf x = - \int_\Omega a  |\nabla \phi_\varepsilon|^2  
d \mathbf x 
- \int_\Omega \nabla a \cdot \phi_\varepsilon  \nabla
\phi_\varepsilon  d \mathbf x, \\
\int_\Omega a \Delta \phi _\varepsilon d \mathbf x =
- \int_\Omega \nabla a \cdot \nabla \phi _\varepsilon  d \mathbf x. 
\end{eqnarray*}
Putting all together,
\begin{eqnarray*}
-  \varepsilon \int_\Omega |\Delta \phi _\varepsilon|^2  d \mathbf x
+ \int_\Omega \left[  {1\over 2} {\rm div}(\mathbf v_l)  -
a \right] |\nabla \phi_\varepsilon|^2  d \mathbf x
+ {1\over 2} \int_{\partial \Omega}  |\nabla \phi_{\varepsilon}|^2 \mathbf v_l 
\cdot \mathbf n d S_{\mathbf x} = \\
\int_\Omega \nabla \! \left[ - {\rm div}(\mathbf v_l)   + a \right]  \cdot
\phi_\varepsilon  \nabla \phi_\varepsilon  d \mathbf x  
 -   \int_\Omega \nabla a \cdot \nabla \phi _\varepsilon  d \mathbf x 
 -   \int_\Omega v_{l,j,x_k} \phi_{\varepsilon, x_j}  \phi_{\varepsilon, x_k}  
 d \mathbf x. 
\end{eqnarray*}
We know that $0\leq \phi _\varepsilon \leq 1$.
Therefore, 
\begin{eqnarray*}
\int_\Omega \left[  -{1\over 2} {\rm div}(\mathbf v_l)  +
a  \right] |\nabla \phi_\varepsilon|^2  d \mathbf x
\leq \left[\|\nabla {\rm div}(\mathbf v_l)\|_{[L^2]^n} 
+ 2  \|\nabla a\|_{[L^2]^n} \right] \| \nabla  \phi_\varepsilon \|_{L^2} \\
+  \int_\Omega |v_{l,j,x_k} \phi_{\varepsilon, x_j}  \phi_{\varepsilon, x_k}| 
d \mathbf x.
\end{eqnarray*}

If $\|\nabla \mathbf v_{l} \|_{[L^\infty]^{n^2}}$ is small enough compared
to $a_{\rm min}$
\begin{eqnarray*}
{1\over 2} a_{\rm min}  \| \nabla \phi_\varepsilon \|_{L^2} \leq 
\|\nabla {\rm div}(\mathbf v_l)\|_{[L^2]^n} + 2  \|\nabla a\|_{[L^2]^n}.
\end{eqnarray*}
As a result, a subsequence $\phi_{\varepsilon'}$ converges weakly in 
$H^1(\Omega)$ to a limit $\phi$, strongly in $L^2(\Omega)$, and 
pointwise in $\Omega$. Traces of $\phi$ are defined on $\partial \Omega$ 
as belonging to $L^2(\partial \Omega)$, and are weak limits of traces
of  $\phi_{\varepsilon'}$. Passing to the limit in the variational formulation
for (\ref{eq_phie}), we see that $\phi \in H^1(\Omega)$ is a solution with 
$\epsilon =0$ which inherits the bounds established.

{\it Uniqueness.} Assume we have two solutions $\phi_1, \phi_2 
\in H^1(\Omega)$ and set $\psi = \phi_1-\phi_2$. Writing down the 
variational equation, substracting, and using $\psi \in H^1(\Omega)$ 
as a test function, we get
\begin{eqnarray*} 
- {1\over 2} \int_{\partial \Omega} (\mathbf v_l \cdot \mathbf n) |\psi|^2 
\, d S_{\mathbf x} +  \int_\Omega \left[ - {1\over 2} {\rm div}(\mathbf v_l)  +
a \right] |\psi|^2  \, d \mathbf x  = 0,
\end{eqnarray*}
which implies $\phi_1=\phi_2$ since $ - {1\over 2} {\rm div}(\mathbf v_l)  +
a>0$ in $\Omega$ and $ -(\mathbf v_l \cdot \mathbf n) \geq 0$. 
$\square$ \\
 
To obtain $L^q$ regularity, we need to adapt Lemma 3.1 from \cite{beirao} to 
our boundary conditions and thin film geometry.

{\bf Lemma 3.3.} {\it Set $q \in [2,\infty)$. Let $\Omega$ be a $C^3$ domain
and $\phi \in W^{3,q}(\Omega)$ such that ${\partial \phi \over \partial
\mathbf n} =0$ on $\partial\Omega$. Assume that $\Omega$ is a thin
domain for which $\mathbf n \sim \mathbf e_n$. Then, for every $\delta >0$,
\begin{eqnarray*}
- \int_{\Omega} \Delta(\nabla \phi) \cdot \left[
(|\nabla \phi|^2 + \delta )^{(q-2)/2} \nabla \phi 
\right] \, d \mathbf x \geq 0. \label{ap1}
\end{eqnarray*}
In particular,
\begin{eqnarray*}
- \int_{\Omega} \Delta(\nabla \phi) \cdot \left[
|\nabla \phi|^{(q-2)} \nabla \phi 
\right] \, d \mathbf x \geq 0. \label{ap2}
\end{eqnarray*}
}
{\bf Proof.} Integrating by parts we find
\begin{eqnarray*}
- \int_{\Omega} \Delta(\nabla \phi) \cdot \left[
(|\nabla \phi|^2 + \delta )^{(q-2)/2} \nabla \phi \right]
\, d \mathbf x = \\
\sum_{i=1}^n \int_{\Omega} {\partial  \over \partial x_i}\left(\nabla \phi \right)
\cdot {\partial   \over \partial   x_i}
\left[  (|\nabla \phi|^2 + \delta )^{(q-2)/2} \nabla \phi \right] \, d \mathbf x \\
- \sum_{i,k=1}^n \int_{\partial \Omega} {\partial^2  \phi \over \partial x_i 
\partial x_k} {\partial \phi \over \partial x_k}  n_i (|\nabla \phi|^2 + \delta)^{(q-2)/2}
\, d S_{\mathbf x}.
\end{eqnarray*}
The first integral is nonnegative since
\begin{eqnarray*}
\sum_{i=1}^n  {\partial  \over \partial x_i}\left(\nabla \phi \right)
\cdot {\partial  \over \partial  x_i}
\left[  (|\nabla \phi|^2 + \delta )^{(q-2)/2} \nabla \phi \right] \, d \mathbf x
= \\ (|\nabla \phi|^2 + \delta )^{(q-2)/2}  \sum_{i,k=1}^n 
({\partial^2  \phi \over \partial x_i \partial x_k})^2
+ {q-2 \over 4} (|\nabla \phi|^2 + \delta )^{(q-4)/2} 
\sum_{i=1}^n ({\partial \over \partial x_i} |\nabla \phi|^2)^2 \geq 0.
\end{eqnarray*}
To estimate the boundary integral, in boundary regions
where $\mathbf n = \mathbf e_n$,
we have ${\partial \phi \over \partial x_n} = 0$, thus
\begin{eqnarray*}
\sum_{k=1}^n {\partial^2  \phi \over \partial x_n \partial x_k}
{\partial \phi \over \partial x_k}    =
\sum_{k=1}^{n-1} {\partial^2  \phi \over \partial x_n \partial x_k}
{\partial \phi \over \partial x_k}   + 
{\partial^2  \phi \over \partial x_n^2} {\partial \phi \over \partial x_n}  =0  
\end{eqnarray*}
since ${\partial^2  \phi \over \partial x_n \partial x_k}
{\partial \phi \over \partial x_k} = 0$ for $k\neq n$.

In fact, we can extend the result to more general
situations using local parametrizations of the boundary. Given
$\mathbf x_0 \in \partial \Omega$, by an orthormal change
of coordinates we may assume that $\mathbf n(\mathbf x_0)$
points in the $x_n$ direction and that the principal directions
of $\partial \Omega$ are parallel to the $x_i$ directions,
$i=1,\ldots,n-1$. The boundary term then vanishes. $\square$ \\

{\bf Theorem 3.4 ($L^q$ regularity). }{\it Under the hypotheses of 
Theorems 3.1 and 3.2, if $\partial \Omega$ has $C^3$ regularity, 
$\Omega$  is a thin domain as in Lemma 3.3, $a \in W^{1,q}(\Omega)$
and ${\rm div}(\mathbf v_l) \in L^\infty(\Omega) \cap W^{1,q}(\Omega)$,
 $n<q<\infty$, then
$\|\nabla \phi \|_{[L^q]^n}$ is bounded from above by a constant
depending on $\|\nabla a\|_{[L^q]^n}$ and 
$\|\nabla {\rm div}(\mathbf v_l)\|_{[L^q]^n}$.  
} 

{\bf Proof.}
By elliptic regularity, 
$\phi _\varepsilon \in W^{3,q}(\Omega)$, since the 
right hand side in (\ref{eq_phie}) belongs to 
$W^{1,q}(\Omega)$. As in \cite{beirao}, we take the derivative 
of both sides of (\ref{eq_phie}) with respect to $x_k$, multiply
by $h(\phi_\varepsilon) \phi_{x_k}$ for $h(\phi_\varepsilon)
= (|\nabla \phi_\varepsilon|^2 + \delta)^{(q-2)/2}$,
add  with respect to $k$ and integrate over $\Omega$,
to find
\begin{eqnarray*}
- \varepsilon \int_{\Omega} \Delta(\nabla \phi_\varepsilon)
\cdot h(\phi_\varepsilon) \nabla \phi_\varepsilon \, d \mathbf x + 
\int_{\Omega} a h(\phi_\varepsilon) |\nabla \phi_\varepsilon|^2 \, 
d \mathbf x  +
\int_{\Omega} \nabla a \cdot h(\phi_\varepsilon) \phi_\varepsilon 
 \nabla \phi_\varepsilon  \, d \mathbf x  - \\
\int_{\Omega} \!\! v_{l,i} \phi_{\varepsilon,x_i x_k} h(\phi_\varepsilon) 
\phi_{\varepsilon, x_k}   d \mathbf x 
\!-\! \int_{\Omega} \!\!  v_{l,i,x_k} \phi_{\varepsilon, x_i} h(\phi_\varepsilon)
\phi_{\varepsilon, x_k}     d \mathbf x  
\!-\! \int_{\Omega} \!\!  {\rm div}(\mathbf v_l) h(\phi_\varepsilon) |\nabla 
\phi_\varepsilon|^2  d \mathbf x  \\ -
\int_{\Omega} \nabla({\rm div}(\mathbf v_l)) \cdot h(\phi_\varepsilon) 
\phi_\varepsilon \nabla \phi_\varepsilon \, d \mathbf x = 
\int_{\Omega} \nabla a \cdot h(\phi_\varepsilon)  
 \nabla \phi_\varepsilon  \, d \mathbf x,
\end{eqnarray*}
where sum over repeated subindices is intended.
By the  Lemma 3.3, the first term is nonnegative in the thin domains
we consider. The fourth term becomes
\begin{eqnarray*}
{1\over q} \int_{\Omega} {\rm div}(\mathbf v_l)(|\nabla \phi_\varepsilon|^2 
+ \delta)^{q/2} d \mathbf x 
- {1\over q} \int_{\partial \Omega} (|\nabla \phi_\varepsilon|^2 
+ \delta)^{q/2}  \mathbf v_{l} \cdot \mathbf n \, dS_{\mathbf x}
\end{eqnarray*}
where the second term is nonnegative. Combing this information we find
\begin{eqnarray*}
 \int_{\Omega} a h(\phi_\varepsilon) |\nabla \phi_\varepsilon|^2 \, 
d \mathbf x 
\leq - {1\over q} \int_{\Omega} {\rm div}(\mathbf v_l)(|\nabla \phi_\varepsilon|^2 
+ \delta)^{q/2} d \mathbf x \\
+ \int_{\Omega}  v_{l,i,x_k} \phi_{\varepsilon, x_i} \phi_{\varepsilon, x_k}  
h(\phi_\varepsilon) \, d \mathbf x
+ \int_{\Omega} {\rm div}(\mathbf v_l) h(\phi_\varepsilon) |\nabla \phi_\varepsilon|^2 
\, d \mathbf x \\
+ \int_{\Omega} \nabla({\rm div}(\mathbf v_l)) \cdot h(\phi_\varepsilon) 
\phi_\varepsilon \nabla \phi_\varepsilon \, d \mathbf x
+ \int_{\Omega} \nabla a \cdot (1-\phi_\varepsilon )h(\phi_\varepsilon)  
 \nabla \phi_\varepsilon  \, d \mathbf x.
\end{eqnarray*}
Letting $\delta \rightarrow 0$
\begin{eqnarray*}
 \int_{\Omega}  a |\nabla \phi_\varepsilon|^q \, d \mathbf x 
\leq (1 - {1\over q}) \int_{\Omega}   {\rm div}(\mathbf v_l) |\nabla 
\phi_\varepsilon|^q d \mathbf x +  
 \int_{\Omega}   v_{l,i,x_k} \phi_{\varepsilon, x_i} \phi_{\varepsilon, x_k}  
 |\nabla \phi_\varepsilon|^{q-2} \, d \mathbf x  \\
+ \int_{\Omega} \nabla({\rm div}(\mathbf v_l)) \cdot  
|\nabla \phi_\varepsilon|^{q-2} \phi_\varepsilon \nabla \phi_\varepsilon 
\, d \mathbf x
+ \int_{\Omega} \nabla a \cdot (1-\phi_\varepsilon ) |\nabla \phi_\varepsilon|^{q-2}
 \nabla \phi_\varepsilon  \, d \mathbf x.
\end{eqnarray*}
For small enough $\| \nabla \mathbf v_l \|_{[L^\infty]^{n^2}}$
\begin{eqnarray*}
{1 \over 2} a_{\rm min}  \int_{\Omega}  |\nabla \phi_\varepsilon|^q \, d \mathbf x 
\leq \left[ \| \nabla({\rm div}(\mathbf v_l)) \|_{L^q} + \| \nabla a \|_{L^q} \right]
        \| |\nabla \phi_\varepsilon| \|_{L^{q}}^{q-1},
\end{eqnarray*}
that is,
\begin{eqnarray*}
{1 \over 2} a_{\rm min }\|  \nabla \phi_\varepsilon   \|_{L^q}
\leq  \| \nabla({\rm div}(\mathbf v_l)) \|_{L^q} + \| \nabla a \|_{L^q}.  
\label{bound_phi2}
\end{eqnarray*}
Letting $\varepsilon \rightarrow 0$, the limit $\phi$ inherits this bound.
$\square$ \\

Once we have established these results, we obtain the following corollaries
for (\ref{vfractionl_st}) and (\ref{vfractionl_stinfty}). \\

{\bf Corollary 3.5.} {\it 
Let $\Omega \subset \mathbb R^n$, $n=2,3,$ be a thin open, bounded 
subset, with $C^3$ boundary $\partial \Omega$. 
Let $\mathbf v_l \in [H^2(\Omega) \cap C(\overline{\Omega})]^n$
such that  ${\rm div}(\mathbf v_l)\leq 0$ in $\Omega$,
${\rm div}(\mathbf v_l) \in L^\infty(\Omega) \cap W^{1,q}(\Omega)$,
$n<q<\infty$, and $\mathbf v_l \cdot   \mathbf n \leq 0$ a.e. on $\partial \Omega$.  
Let $c \in W^{1,q}(\Omega) \cap C(\overline{\Omega})$  be a strictly 
positive function. We assume that $c$ is bounded from below by a positive 
constant $c_{\rm min}$, that is, $c \geq c_{\rm min}>0$ in $\Omega$.
We assume that $ \nabla \mathbf v_l \in [L^{\infty}(\Omega)]^{n^2}$ with 
$\|\nabla \mathbf v_l\|_{[L^{\infty}]^{n^2}}$ small enough compared to 
$c_{\rm min}$.
Then, given positive constants $k_b$ and  $K_b$,  there exists 
a  solution $\phi_l \in L^2(\Omega)$ of  (\ref{vfractionl_st}) in the 
sense of distributions. Moreover,
\begin{itemize}
\item $0 \leq \phi_l \leq 1$ on $\Omega$, and $\phi$ does not vanish in sets
of positive measure.
\item $\phi_l \in H^1(\Omega)$ is the unique solution of the
variational formulation in $H^1(\Omega)$ and 
\begin{eqnarray*}
{1\over 2} {c_{\rm min} k_b \over \| c\|_{L^\infty} + K_b}
 \| \nabla \phi  \|_{L^2} \leq 
\|\nabla {\rm div}(\mathbf v_l)\|_{[L^2]^n} + 2 k_b K_b \|\nabla c\|_{[L^2]^n}.
\end{eqnarray*}
\item $\nabla \phi_l \in L^q(\Omega)$ and
\begin{eqnarray*}
{1\over 2} {c_{\rm min} k_b \over \| c\|_{L^\infty} + K_b}
 \| \nabla \phi  \|_{L^q} \leq 
\|\nabla {\rm div}(\mathbf v_l)\|_{[L^q]^n} + 2 k_b K_b \|\nabla c\|_{[L^q]^n}.
\end{eqnarray*}
\end{itemize}
}

{\bf Corollary 3.6.} {\it 
Let $\Omega \subset \mathbb R^n$, $n=2,3,$ be a thin open, bounded 
subset, with $C^3$ boundary $\partial \Omega$. 
Let $\mathbf v_l \in [H^2(\Omega) \cap C(\overline{\Omega})]^n$
such that  ${\rm div}(\mathbf v_l)\leq 0$ in $\Omega$,
${\rm div}(\mathbf v_l) \in L^\infty(\Omega) \cap W^{1,q}(\Omega)$,
$n<q<\infty$, and $\mathbf v_l \cdot   \mathbf n \leq 0$ a.e. on $\partial \Omega$.  
We assume that $ \nabla \mathbf v_l \in [L^{\infty}(\Omega)]^{n^2}$ with 
$\|\nabla \mathbf v_l\|_{[L^{\infty}]^{n^2}}$ small enough compared to 
$k_b g_\infty$.
Then, given positive constants $k_b$ and  $g_\infty$,  there exists 
a  solution $\phi_l \in L^2(\Omega)$ of  (\ref{vfractionl_stinfty}) in the 
sense of distributions. Moreover,
\begin{itemize}
\item $0 \leq \phi_l \leq 1$ on $\Omega$ and $\phi$ does not vanish in 
sets of positive measure.
\item $\phi_l \in H^1(\Omega)$ is the unique solution of the
variational formulation in $H^1(\Omega)$ and 
\begin{eqnarray*}
{1\over 2} k_b g_\infty
 \| \nabla \phi  \|_{L^2} \leq 
\|\nabla {\rm div}(\mathbf v_l)\|_{[L^2]^n}.
\end{eqnarray*}
\item $\nabla \phi_l \in L^q(\Omega)$ and
\begin{eqnarray*}
{1\over 2} k_b g_\infty
 \| \nabla \phi  \|_{L^q} \leq 
\|\nabla {\rm div}(\mathbf v_l)\|_{[L^q]^n}.
\end{eqnarray*}
\end{itemize}
}

{\bf Remark.} Once we know $\phi_l$ we can construct 
$\phi_b = 1-\phi_l$. Similar arguments prove the existence of solutions 
to the quasi-stationary version of (\ref{conservationb})
assuming ${\rm div}(\mathbf v_b) \geq 0$, and
$ \mathbf v_b \cdot \mathbf n \geq 0$, and given data $\phi_b=\phi \in (0,1)$ 
on $\partial \Omega$. Then, we would set $\phi_l = 1- \phi_b$.

\section{The Stokes problem}
\label{sec:stokes}

Consider the vector system
\begin{eqnarray}
\mu_b \Delta \mathbf v_b +  {\mu_b \over 3}  \nabla {\rm div} 
(\mathbf v_b) = \nabla (\Pi \phi_b + p), &  \mathbf x \in \Omega,
 \label{velocity_st} \\
 {\rm div}(\xi(\phi_b) \nabla p) = {\rm div}(\mathbf v_b),
&   \mathbf x \in \Omega,
\label{pressure_st}
\end{eqnarray}
where $\mu_b$ and $\Pi$ are positive constants and $\xi(\phi_b)$ is a
known scalar function. Typically,  $\xi(\phi_b) ={ (1-\phi_b)^2 \over \mu_l} 
\zeta_\infty$, with $\zeta_{\infty}>0$, constant. 
System (\ref{velocity_st})-(\ref{pressure_st}) is a variant of compressible
Stokes equations with additional regularity on the pressure term enforced 
by the coupling.

On the boundary we impose the boundary conditions
\begin{eqnarray} \begin{array}{ll}
[ \mu_b (\nabla \mathbf v_b + \nabla \mathbf v_b^t)
- (2 {\mu_b\over 3} {\rm div}(\mathbf v_b) + ( p + \Pi \phi_b)) \mathbf I] 
\cdot \mathbf n = \mathbf t_{\rm ext}, & 
 \mathbf x \in \partial \Omega^+,  \\[1ex]
\mathbf v_b = 0, &   \mathbf x \in \partial \Omega^-,   \\[1ex]
p = p_{\rm ext} - \pi_{\rm ext}, & \mathbf x \in \partial \Omega,
\end{array} \label{bc_vp_st}
\end{eqnarray}
$\mathbf n$ being the unit outer normal.  \\


{\bf Theorem 4.1 (Existence and uniqueness).} 
{\it Let $\Omega \subset \mathbb R^n$, 
$n=2,3,$ be an open bounded domain with $C^1$ boundary $\partial 
\Omega$.  Consider the solution $\phi_b \in L^2(\Omega)$ of
(\ref{vfractionl_gen}) provided by Theorem 3.1. Assume that
$ \mathbf t_{\rm ext} \in [L^2(\partial \Omega)]^n$ and $p_{\rm ext} - 
\pi_{\rm ext} \in H^{1/2}(\partial \Omega) \cap L^2(\partial \Omega)$.
Given  positive constants $\mu_b$, $\Pi$ and a strictly positive 
function $\xi(\phi_b) \in C^1([0,1])$ there exists a unique solution 
$(\mathbf v_b,  p) \in [H^1(\Omega)]^n \times H^1(\Omega)$ of 
(\ref{velocity_st})-(\ref{pressure_st}) with boundary conditions 
(\ref{bc_vp_st}).
Moreover, 
 $\|\mathbf v_b\|_{[H^1]^n}$ and $\|p \|_{H^1}$ are bounded from
above by constants depending on $\mu_b$, $\xi$, $\Pi$, $\Omega$,
$\|\phi_b\|_{L^2}$, $\|\mathbf t_{\rm ext}\|_{[L^2(\partial \Omega)]^n}$, 
$\| p_0\|_{H^1}$, $p_0$ being a $H^1(\Omega)$ extension of  
$p_{\rm ext} - \pi_{\rm ext}.$ 
}

{\bf Proof.} {\it Existence and uniqueness.}
By the theory of traces \cite{gagliardo,lions},
any $p_{\rm ext} - \pi_{\rm ext} \in H^{1/2}(\partial \Omega)$ is the trace
of a function $p_0\in H^1(\Omega).$ Setting $p = \tilde p + p_0$,
\begin{eqnarray} \begin{array}{ll}
 {\rm div}(\xi(\phi_b) \nabla \tilde p) = {\rm div}(\mathbf v_b) -  
 {\rm div}(\xi(\phi_b) \nabla p_0),   &  \mathbf x \in \Omega, 
  \\[1ex]
 \tilde p = 0, & \mathbf x \in \partial \Omega.
\end{array}  \label{pressure_tilde} 
\end{eqnarray}
The function $\phi_b$ is defined on $\partial \Omega$ as a 
$H^{-1/2}(\partial \Omega)$ function \cite{lions}.

We define the bilinear forms
\begin{eqnarray*}
b_1(\mathbf v, \mathbf w) =  \int_{\Omega}
[ \mu_b \nabla \mathbf v \nabla \mathbf w^t  + {\mu_b \over 3} 
{\rm div}(\mathbf v) {\rm div} (\mathbf w ) ] \, d \mathbf x, \quad
b_2(\tilde p, \mathbf w) = -
\int_{\Omega}   \tilde p  \, {\rm div}(\mathbf w) \, d \mathbf x, \\
c_1(\tilde p, q) = \int_{\Omega} \xi(\phi_b) \nabla \tilde p \cdot 
\nabla q \, d \mathbf x, \quad
c_2(\mathbf v, q) = \int_{\Omega} {\rm div}(\mathbf v) q \, d \mathbf x,
\end{eqnarray*}
on $H= [H^1_{0,-}(\Omega)]^n \times H^1_0(\Omega)$. Here, $H^1_0$
is the standard Sobolev space of $H^1$ functions vanishing of 
$\partial \Omega$ and $H^1_{0,-}$ the space of $H^1$ functions 
vanishing only on $\partial \Omega^-$. For $[\mathbf v,\tilde p], 
[\mathbf w,q] \in H$
\begin{eqnarray*}
a([\mathbf v,\tilde p],[\mathbf w,q]) = b_1(\mathbf v, \mathbf w)  
+ b_2(\tilde p, \mathbf w) + c_1(\tilde p, q) + c_2(\mathbf v, q) 
\label{st_biblinear}
\end{eqnarray*}
defines a bilinear form, continuous on $H \times H$.
We multiply the equations by test functions $[\mathbf w,q]  \in H$, 
integrate by parts and add to get the weak formulation
\begin{eqnarray*}
a([\mathbf v,\tilde p],[\mathbf w,q]) = L([\mathbf w,q]), \label{st_variational}
\end{eqnarray*}
of system (\ref{velocity_st}), (\ref{pressure_tilde}), (\ref{bc_vp_st}), 
 where  $L([\mathbf w,q])$ is given by
\begin{eqnarray*}
 L([\mathbf w,q]) \!=\! \int_{\Omega} \hskip -1mm 
 ( \Pi \phi_b \!+\! p_0)\,{\rm div} (\mathbf w) \, d \mathbf x  \!+\! 
 \int_{\partial \Omega}  \hskip -1mm \mathbf t_{\rm ext}  
 \cdot \mathbf w \, d S_{\mathbf x} \!-\! 
 \int_{\Omega} \hskip -1mm \xi(\phi_b) \nabla p_0\cdot  \nabla q \, d \mathbf x. 
 \label{st_linear}
\end{eqnarray*}
Notice that we have
$_{H^{-1}}< {\rm div}(\xi(\phi_b) \nabla p_0), q >_{H_0^1}
=  - \int_{\Omega} \xi(\phi_b) \nabla p_0 \cdot \nabla q \, d \mathbf x
+ \int_{\partial \Omega} \xi(\phi_b) q \nabla p_0 \cdot \mathbf n \, 
d S_{\mathbf x}$
where the boundary term vanishes because $q \in H^1_0(\Omega)$.
Similarly, $ <\Pi \nabla \phi_b, \mathbf w > = 
-  \int_{\Omega}  \Pi \phi_b \,{\rm div} (\mathbf w) \, d \mathbf x
+ < \Pi \phi_b \mathbf n, \mathbf w >$ for $\mathbf w \in [H^1_{0,-}(\Omega)]^n$,
see \cite{lions}.

Since $\phi_b \in L^2(\Omega)$, $p_0 \in L^2(\Omega)$, 
$ \mathbf t_{\rm ext}  \in [L^2(\partial \Omega)]^n$ and 
$ \xi(\phi_b) \nabla p_0  \in L^2(\Omega),$ $L$ defines a continuous 
linear form on $H$. Setting 
$[\mathbf w,q]=[\mathbf v,\tilde p]$ and $\xi_{\rm min} = {\rm min}|\xi|$, 
we see that two terms cancel and 
\begin{eqnarray*}
a([\mathbf v,\tilde p],[\mathbf v,\tilde p]) = b_1(\mathbf v, \mathbf v)  +
c_1(\tilde p, \tilde p)  \geq 
\mu_b \int_{\Omega} |\nabla \mathbf v|^2 \, d \mathbf x   + \xi_{\rm min} 
\int_{\Omega} |\nabla  \tilde p |^2 \, d \mathbf x  \\ \geq  
\mu_b P_{\partial \Omega^-} \|\mathbf v\|^2_{[H^1]^n} + \xi_{\rm min} 
 P_{\partial \Omega} \|\tilde p \|^2_{H^1}
\end{eqnarray*}
thanks to Poincar\'e inequalities for functions vanishing on parts on 
the boundary with positive measure \cite{brezis,raviart2}. 
$P_{\partial \Omega}>0$ and  $P_{\partial \Omega^-}>0$ denote the 
corresponding Poincar\'e constants for $\Omega$. Lax-Milgram 
theorem ensures the  existence of a unique solution $(\mathbf v_b, \tilde p)$ 
in $H$ \cite{brezis}.

{\it Bounds.}
Taking  $[\mathbf w,q]=[\mathbf v,\tilde p]$ in the variational equation 
and using coercivity, we see that
\begin{eqnarray*}
\mu_b P_{\partial \Omega^-} \|\mathbf v\|^2_{[H^1]^n} + 
\xi_{\rm min} P_{\partial \Omega} \|\tilde p \|^2_{H^1} \leq 
[ \Pi \| \phi_b\|_{L^2} + \| p_0\|_{L^2}]    \|\mathbf v\|_{[H^1]^n}   \\
+ T_\Omega \|\mathbf t_{\rm ext}  
\|_{[L^2(\partial \Omega)]^n} \|\mathbf v\|_{[H^1]^n}
+ \| \xi(\phi_b) \nabla p_0  \|_{L^2} \|\tilde p \|_{H^1}.
\end{eqnarray*}
Setting $\xi_{\rm max} = {\rm max}|\xi|$, this implies 
\begin{eqnarray}
{\rm min}(\mu_b P_{\partial \Omega^-}, \xi_{\rm min} P_{\partial \Omega}
) \left[ \|\mathbf v\|_{[H^1]^n} +  \|\tilde p \|_{H^1} \right]
\leq \nonumber \\ 
2[ \Pi \| \phi_b\|_{L^2}+ T_\Omega \|\mathbf t_{\rm ext} 
\|_{[L^2(\partial \Omega)]^n}  + (\xi_{\rm max}+1)
 \| p_0\|_{H^1}].    \label{bound_vp1}
\end{eqnarray}
$\square$

In practice, we will replace $\xi(\phi_b)$ by $\xi_\infty={ \phi_\infty^2 
\over \mu_l} \zeta_\infty >0$, with $\phi_l = \phi_\infty \in (0,1)$, to simplify 
the obtention of higher order bounds. In fact, simulations and asymptotic 
studies suggest that $\phi_l$ remains close to a constant.
Then  equation (\ref{pressure_st}) becomes
\begin{eqnarray}
 {\rm div}(\xi_\infty \nabla p) = \xi_\infty \Delta p = {\rm div}(\mathbf v_b),
&   \mathbf x \in \Omega.
\label{pressure_stinfty}
\end{eqnarray}

{\bf Theorem 4.2 (Regularity).} {\it
Let $\Omega \subset \mathbb R^n$, $n=2,3,$ be an open bounded 
domain with  $C^2$ boundary $\partial \Omega$.  
Consider the solution $\phi_b \in L^\infty(\Omega)$ of
(\ref{vfractionl_gen}) provided by Theorem 3.1.  Assume that
$ \mathbf t_{\rm ext} \in [L^2(\partial \Omega)\cap L^q(\partial \Omega)]^n$ 
and $p_{\rm ext} - \pi_{\rm ext} \in H^{3/2}(\partial \Omega) \cap 
W^{2-1/q,q}(\partial \Omega)$, $n<q<\infty$.
Given  positive constants $\mu_b$, $\Pi$, $\xi_\infty,$ 
there exists a unique solution 
$(\mathbf v_b,  p) \in [H^1(\Omega)]^n \times H^2(\Omega)$  of 
(\ref{velocity_st}), (\ref{pressure_stinfty}) with 
boundary conditions (\ref{bc_vp_st}). Moreover, 
\begin{itemize}
\item $\|\mathbf v_b\|_{[H^1]^n}$ and $\|p \|_{H^1}$ are bounded from
above by constants depending on $\mu_b$, $\xi_\infty$, $\Pi$, $\Omega$,
$\|\phi_b\|_{L^2}$, $\|\mathbf t_{\rm ext} \|_{[L^2(\partial \Omega)]^n}$, 
$\|p_0\|_{H^1}$. Additionally, $\|p \|_{H^2}$ is bounded
by constants depending on $\xi_\infty$,  $\Omega$,
$\| {\rm div}(\mathbf v_l) \|_2$ and $\|p_0\|_{H^2}$, 
$p_0$ being a $H^2\cap W^{2,q}$ extension of  
$p_{\rm ext} - \pi_{\rm ext}.$
\item ${\rm div}(\mathbf v_b) \in L^q(\Omega)$, $\mathbf v_b \in 
[W^{1,q}(\Omega)]^n$ and $p \in W^{2,q}(\Omega)$. The norms
$\|\mathbf v_b \|_{[W^{1,q}]^n}$ and $\|p \|_{W^{2,q}}$
are bounded from above by constants depending
on  $\mu_b$, $\xi_\infty$, $\Pi$, $\Omega$,
$\|\phi_b\|_{L^q}$,  $\|\mathbf t_{\rm ext}\|_{[L^q(\partial \Omega)]^n}$
and $\|p_0\|_{W^{2,q}}$. Moreover, 
$\mathbf v_b \in C(\overline{\Omega})$ and $p \in C^1(\overline{\Omega})$.
\end{itemize}
}

{\bf Proof.}
Any $p_{\rm ext} - \pi_{\rm ext} \in H^{3/2}(\partial \Omega) 
\cap W^{2-1/q,q}(\partial \Omega)$ is the trace of a function 
$p_0 \in H^2(\Omega) \cap W^{2,q}(\partial \Omega)$. 
Existence of a unique solution 
$(\mathbf v_b,  p) \in [H^1(\Omega)]^n \times H^1(\Omega)$ follows as 
in the proof of Theorem 4.1. We obtain estimate (\ref{bound_vp1}) replacing 
$\xi_{\rm min}$ and $\xi_{\rm max}$ by $\xi_\infty$. 

We set again $p=\tilde p +p_0$.
Elliptic regularity for the laplacian $\xi_\infty \Delta \tilde p = {\rm div} 
(\mathbf v_b) - \xi_\infty \Delta p_0 \in L^2(\Omega)$  provides $H^2(\Omega)$ 
regularity  and estimates for $p$. By Sobolev injections \cite{adams,brezis}, 
$p \in L^\infty(\Omega) \cap C(\overline{\Omega})$ and 
$\nabla p \in L^q(\Omega)$, $q<q^*$, $q^*=\infty$ if $n=2$, or $q^*= 6$
 if $n=3$.

To increase regularity, notice that we can rewrite (\ref{velocity_st}) and 
(\ref{pressure_stinfty}) as an elliptic system in $\Omega$ \cite{adn2,kozlov}:
\begin{eqnarray} \begin{array}{lcl}
\mu_b \Delta v_{b,1} + {\mu_b\over 3} {\partial^2 v_{b,1} \over \partial x_1^2} 
+ \ldots + {\mu_b\over 3} {\partial^2 v_{b,1} \over \partial x_1\partial x_n} 
-{\partial \tilde p \over \partial x_1} & = & \Pi {\partial \phi_b \over \partial x_1}
+ {\partial p_0 \over \partial x_1}, \\
\ldots & = & \ldots, \\
\mu_b \Delta v_{b,n} + {\mu_b\over 3} {\partial^2 v_{b,1} \over \partial x_n\partial x_1} 
+\ldots + {\mu_b\over 3} {\partial^2 v_{b,n} \over \partial x_n^2} 
-{\partial \tilde p \over \partial x_n} & = & \Pi {\partial \phi_b \over \partial x_n}
+ {\partial p_0 \over \partial x_n},\\
\xi_\infty \Delta \tilde p - {\partial v_{b,1} \over \partial x_1} - \ldots
- {\partial v_{b,n} \over \partial x_n}  & = & - \xi_\infty \Delta p_0.
\end{array} \label{system_vp} \end{eqnarray}
The right hand side belongs to $ [W^{-1,q}(\Omega)]^{n+1}, $ $n< q < \infty$.
Elliptic regularity for general boundary conditions \cite{adn2,kozlov} implies 
then that $v_{b,j} \in W^{1,q}(\Omega)$, $j=1, \ldots, n$, $n<q <
\infty$. Therefore, $v_{b,j} \in C(\overline{\Omega})$, $j=1, \ldots, n$. Moreover,
\begin{eqnarray} \begin{array}{r}
\|\mathbf v_b \|_{[W^{1,q}]^n} + \|\tilde p \|_{W^{1,q}} \leq  \\ [1ex]
K_q(\Omega,\mu_b,\xi_\infty) \left[ \Pi \|\phi_b\|_{L^q} +
\|p_0\|_{W^{1,q}} + \|\mathbf t_{\rm ext}\|_{[L^q(\partial \Omega)]^n} \right].
\end{array} \label{bound_vp2}
\end{eqnarray}

Now, elliptic regularity for the laplacian (\ref{pressure_stinfty}) with a right hand 
side and boundary data in $L^q$ imply $p \in W^{2,q}(\Omega)$  and 
$\nabla  p \in W^{1,q}(\Omega)$ with
\begin{eqnarray}
 \|  p \|_{W^{2,q}} \leq K_{2q}(\Omega) [ {1 \over \xi_\infty} \|\mathbf v_b \|_{[W^{1,q}]^n}
  + \| p_{\rm ext} -\pi_{\rm ext} \|_{L^q(\partial \Omega)} ].
\label{bound_p2}
\end{eqnarray}
In particular, $\nabla p  \in C(\overline{\Omega})$ by Sobolev injections
since $q>n$. 
$\square$
\\

Even if we take $q=\infty$ here, we do not get
$W^{1,\infty}$ regularity for $\mathbf v_b$ assuming just
$\phi_b \in L^\infty$, because elliptic regularity theory does not
ensure that fact. We only have
${\rm div(\mathbf v_b)} \in L^q(\Omega)$, $1\leq  q < \infty$.
Next, we give conditions for $W^{1,\infty}$ regularity. \\

{\bf Theorem 4.3 (Higher Regularity).} {\it Under the hypotheses of Theorem
4.2, let us assume that $\phi_b \in H^1(\Omega)$. Then
$(\mathbf v_b,  p) \in [H^2(\Omega)]^n \times H^2(\Omega)$, and
$\|\mathbf v_b\|_{[H^2]^n}$ and  $\|p\|_{H^2}$ are bounded from
above by constants depending on $\mu_b$, $\xi_\infty$, $\Pi$, $\Omega$,
$\|\phi_b\|_{H^1}$, $\|\mathbf t_{\rm ext} \|_{[L^2(\partial \Omega)]^n}$
and $\|p_0\|_{H^2}$.

If $\Omega$ has $C^3$ regularity, 
$p_{\rm ext} - \pi_{\rm ext} \in  H^{5/2}(\partial \Omega)$ and $p_0$ is 
a $H^3 (\Omega)$ extension of  $p_{\rm ext} - \pi_{\rm ext},$  
then $p \in H^3(\Omega)$ and $\|p  \|_{H^3}$  
is bounded from above by a constant depending on $\xi_\infty$, $\Omega$,
$\|\mathbf v_b\|_{[H^2]^n}$ and $\|p_0  \|_{H^3}$.

Furthermore, if $\phi_b \in W^{1,q}(\Omega)$, $n<q<\infty$, then
$(\mathbf v_b,  p) \in [W^{2,q}(\Omega)]^n \times W^{2,q}(\Omega)$, and
$\|\mathbf v_b\|_{[W^{2,q}]^n}$ and $\|p\|_{W^{2,q}}$ are bounded from
above by constants depending on $\mu_b$, $\xi_\infty$, $\Pi$, $\Omega$,
$\|\phi_b\|_{W^{1,q}}$, $\|\mathbf t_{\rm ext} \|_{[L^q(\partial \Omega)]^n}$
and $\|p_0\|_{W^{2,q}}$.

If $p_{\rm ext} - \pi_{\rm ext} \in  W^{3-1/q,q}(\partial \Omega)$ and $p_0$ is 
a $W^{3,q}(\Omega)$ extension of  $p_{\rm ext} - \pi_{\rm ext},$  
then $p \in W^{3,q}(\Omega)$ and $\|p \|_{W^{3,q}}$  
is bounded from above by a constant depending on $\xi_\infty$, $\Omega$,
 $\|\mathbf v_b\|_{[W^{3,q}]^n}$ and $\|p_0 \|_{W^{3,q}}$.
} 

{\bf Proof.}
These claims are straightforward consequences of Theorem 4.2
and elliptic regularity theory  \cite{adn,gilbart} applied first to (\ref{system_vp})
and  then to (\ref{pressure_stinfty}).
The $H^2$ and $H^3$ norms of the solutions satisfy
\begin{eqnarray*} \begin{array}{l}
\| \mathbf v_b \|_{[H^2]^n} + \| \tilde p \|_{H^2} \leq  K_2(\Omega,\mu_b,\xi_\infty) 
\left[\Pi \| \phi_b \|_{H^1} + \xi_\infty \| p_0 \|_{H^2}+
\|\mathbf t_{\rm ext}\|_{[L^2(\partial \Omega)]^n}  \right], \\[1.5ex]
\|\tilde p\|_{H^3} \leq K_3(\Omega)  \big[{1\over \xi_\infty} 
\| \mathbf v_b \|_{[H^2]^n} + \| p_0 \|_{H^3} \big]. \end{array} 
\end{eqnarray*}

When $\phi_b \in W^{1,q}$, elliptic regularity yields
\begin{eqnarray} \begin{array}{ll}
\|\mathbf v_b \|_{[W^{2,q}]^n} + \|\tilde p \|_{W^{2,q}} \leq & K_{2q}(\Omega,\mu_b,\xi_\infty)
\big[ \Pi \| \nabla \phi_b\|_{L^q} + \|p_0\|_{W^{2,q}} + \\[1ex]
& \|\mathbf t_{\rm ext}\|_{[L^q(\partial \Omega)]^n} \big].
\end{array} \label{bound_vp4}  \\
\begin{array}{ll} 
 \|\tilde p \|_{W^{3,q}} \leq & K_{3q}(\Omega) \big[ {1 \over \xi_\infty} \|\mathbf v_b \|_{[W^{2,q}]^n}  + \| p_0 \|_{W^{3,q}} \big].
\end{array}  \label{bound_p4}
\end{eqnarray}


\section{The convection-reaction-diffusion problem}
\label{sec:crd}

Consider the scalar problem for $c(\mathbf x)$
\begin{eqnarray}
\begin{array}{rcll}
-  d \Delta c + {\rm div}  (\mathbf v_l c)  
   &=&  -  k_c \phi_b { c \over c + K_c}, &
  \mathbf x \in \Omega, \\
  c &=& c_0, &  \mathbf x \in \partial \Omega^-, \\
  {\partial c \over \partial \mathbf n} &=& 0, &  
  \mathbf x \in \partial \Omega^+, \\
\end{array} \label{nutrient_st} 
\end{eqnarray}
where $k_c, K_c, d$ are positive constants,  $\mathbf v_l$
is a known vector function and $\phi_b$ a known scalar function.
Previous asymptotic and computational studies of the biofilm model, 
comparing computational  results to  experimental observations 
\cite{entropy,seminara} indicate that the value of $d$ is very large, 
while $\mathbf v_l$ is fairly small, once dimensions have been 
removed. We will make that assumption here.
\\

{\bf Theorem 5.1.} {\it Let $\Omega \subset \mathbb R^n$, 
$n=2,3$, be a bounded domain with  $C^2$ boundary. Given 
positive constants $k_c, K_c, d, c_0$, a vector function 
$\mathbf v_l \in [W^{1,q}(\Omega)]^n \cap C(\overline{\Omega})$,
$n<q<\infty,$ and a positive function $\phi_b \in  L^\infty(\Omega)$
(not identically zero in a set of positive measure), there 
exists a unique nonnegative solution $c \in H^1(\Omega)$
of (\ref{nutrient_st})  provided $d$ is sufficiently large and 
$\| \mathbf v_l\|_{[W^{1,q}]^n}$ sufficiently small, depending 
on $\Omega$. Moreover,
\begin{itemize}
\item  $\| c \|_{H^1}$ is bounded from above
by constants depending on $d$, $k_c$, $K_c$,
$\Omega$, $c_0$, $\|\mathbf v_l \|_{[L^2]^n}$, 
$\| \phi_b \|_{L^2}.$ 
The same holds for $\| c \|_{H^2}$ and $\| c \|_{L^\infty}=c_{\rm max}$.
\item If ${\rm div}(\mathbf v_l)\leq 0$, then 
$c$ is stricly positive in $\Omega$. 
\item If ${\rm div}(\mathbf v_l)\leq 0$, $0\leq \phi_b \leq 1$ and 
$M= {\rm max} \{ x_n \, | \, (x_1,\ldots,x_{n-1},x_n) \in \Omega\} < 
{2d K_c \over k_c c_{\rm max}} c_0$, then 
$c \geq - {k_c c_{\rm max}\over 2 d K_c} M + c_0 >0$
 in $\overline{\Omega}$.
\end{itemize}}

{\bf Proof.} {\it Existence. }
Starting from $c^{(0)} = c_0$, we consider the iterative 
scheme
\begin{eqnarray*}
\begin{array}{rcll}
-  d \Delta c^{(m+1)} + {\rm div}  (\mathbf v_l c^{(m+1)})  
   &=&  -  k_c { \phi_b  \over c^{(m)} + K_c} c^{(m+1)}, &
  \mathbf x \in \Omega, \\
  c^{(m+1)} &=& c_0, &  \mathbf x \in \partial \Omega^-, \\
  {\partial c^{(m+1)}  \over \partial \mathbf n} &=& 0, &  
  \mathbf x \in \partial \Omega^+, \\
\end{array} 
\end{eqnarray*}
for $m\geq 0$. Setting $c^{(m+1)} = \tilde c^{(m+1)} + c_0$,
these problems admit the variational formulation: Find
$\tilde c^{(m+1)} \in H^1_{0,-}(\Omega)$
such that
\begin{eqnarray*}
b_m(\tilde c^{(m+1)},w)= d \int_{\Omega}  \hskip -1mm \nabla \tilde 
c^{(m+1)} \!\cdot\! \nabla w \, d \mathbf x 
 -\! \int_{\Omega} \hskip -1mm \mathbf v_l \tilde c^{(m+1)} \!\cdot \!
\nabla w \, d \mathbf x 
 +\! \int_{\partial \Omega^+}  \hskip -4mm \tilde c^{(m+1)} w \mathbf v_l 
\!\cdot\! \mathbf n \, dS_{\mathbf x}  \\
 + \int_{\Omega} a_m \tilde c^{(m+1)} w \, d \mathbf x    
 =  c_0 \int_{\Omega} \mathbf v_l \cdot  \nabla w \, d \mathbf x   
 - c_0\int_{\Omega} a_m w \, 
d \mathbf x = L(w),
\end{eqnarray*}
for all $ w \in H^1_{0,-}(\Omega)$, where
$H^1_{0,-}(\Omega) = \{ w \in H^1(\Omega) 
\, |  \, w = 0 \, \mbox{\rm on } \partial \Omega^- \}  $
and $a_m =  k_c  { \phi_b  \over  c^{(m)}  + K_c}$. Notice that
Poincar\'e's inequality holds in $H^1_{0,-}(\Omega)$
\cite{raviart2}. Moreover, $a_m \geq 0$ and $a_m \in L^\infty(\Omega)$
provided $c^{(m)} \geq 0$. In fact, we have the uniform bound
$\| a_m \|_\infty \leq k_c \|\phi_b\|_\infty /K_c$. 

Let us assume that $c^{(m)} \geq 0$ and study the bilinear form 
$b(\tilde  c,w)$ in $H= H^1_{0,-}(\Omega)$. 
Setting $w = \tilde c$, by Sobolev embeddings and trace
inequalities \cite{brezis} we have
\begin{eqnarray*}
\big | \int_{\Omega} \mathbf v_l \cdot \tilde  c \nabla \tilde c
 \, d  \mathbf x \big | \leq  \| \mathbf v_l  
 \|_{[L^\infty]^n} \|\tilde  c \|_{L^2} \|\nabla \tilde  c\|_{[L^2]^n}
 \leq  S_{q,\infty}(\Omega) \| \mathbf v_l  \|_{[W^{1,q}]^n}  
 \| \tilde c\|_{H^1}^2, \\
 \big | \int_{\partial \Omega} |\tilde c|^2 \mathbf v_l \cdot
\mathbf n \, dS_\mathbf x   \big |
 \leq \| \mathbf v_l \|_{[L^\infty]^n} \|\tilde c\|_{L^2(\partial \Omega)}^2
  \leq S_{q,\infty}(\Omega) T_\Omega
   \| \mathbf v_l  \|_{[W^{1,q}]^n} \| \tilde c\|_{H^1}^2,
\end{eqnarray*}
where $T_\Omega, S_{q,\infty}(\Omega)$ denote
trace and embeeding constants. We set $k_{\mathbf v_l, \Omega}
= S_{q,\infty}(\Omega) (1+ T_\Omega) \| \mathbf v_l  \|_{[W^{1,q}]^n}$.
Since $q>n$ and $\Omega$ is bounded, $\mathbf v_l  \in [H^1]^n$.
Using Poincar\'e's inequality in $H$ \cite{raviart2} with Poincare's
constant $P_{\partial \Omega^-} $ and $a_m \geq 0$,  
\begin{eqnarray*}
b_m(\tilde  c, \tilde c) \geq \left[ {d\over 2} \!-\! k_{\mathbf v_l, \Omega} \right]
 \|\nabla \tilde c\|_{[L^2]^n}^2 +
\left[ {P_{\partial \Omega^-}\over 2}  \!-\! k_{\mathbf v_l, \Omega}\right]
\|\tilde  c \|_{L^2}^2
\geq {1\over 4}{\rm min}(d,P_{\partial \Omega^-} ) \|\tilde  c \|_{H}^2,
\end{eqnarray*}
provided $\mathbf v_l$ is small enough to ensure
$P_{\partial \Omega^-}/4  > k_{\mathbf v_l, \Omega} $ and
$d/4 > k_{\mathbf v_l, \Omega}$. 
Thus, $b$ is coercive in $H\times H$, and continuous. 
Moreover, the right  hand side $L$ defines a continuous linear form 
on $H$. Then, Lax Milgram's theorem  \cite{brezis} guarantees the 
existence of a unique solution $\tilde c^{(m+1)} \in H$.

Furthermore, setting $w=\tilde c^{(m)}$ for $m \geq 1$ implies
\begin{eqnarray*}
{1\over 4}{\rm min}(d,P_{\partial \Omega^-} ) \|\tilde c^{(m)}\|_H 
\leq c_0  \,  \left[ \| \mathbf v_l\|_{[L^2]^n} +
{k_c \over K_c} \|\phi_b\|_{L^2}\right].
\end{eqnarray*}
Thus, we can extract a subsequence $\tilde c^{(m_j)}$ 
converging to a limit $\tilde c$ weakly in $H$, strongly in 
$L^2(\Omega)$, and almost everywhere in $\Omega$,
with $\tilde c^{(m_j)} |_{\partial \Omega}$ converging
weakly in $L^2(\partial \Omega)$. We can pass 
to the limit in the variational formulation for $\tilde c^{(m_j)}$ 
obtaining $c = \tilde c + c_0$, with $\tilde c \in H$ a solution of
\begin{eqnarray*}
d \int_{\Omega}  \nabla \tilde c \cdot \nabla w \, d \mathbf x 
- \int_{\Omega} \mathbf v_l \cdot \tilde c \nabla w  \, d \mathbf x 
+ \int_{\partial \Omega^+} \tilde c w \, \mathbf v_l \cdot
\mathbf n \, dS_\mathbf x  + \\
\int_{\Omega} k_c  { \phi_b  \over  c + K_c} \tilde c w \, d \mathbf x    
= c_0 \int_{\Omega} \mathbf v_l \cdot  \nabla w \, d \mathbf x   
- c_0\int_{\Omega} k_c  { \phi_b  \over  c + K_c} w \, d \mathbf x 
\end{eqnarray*}
for $w \in H$.
Convergence of $a_{m_j}\tilde c^{(m_j+1)} w$ to 
$k_c  { \phi_b  \over  c + K_c} \tilde c w$ and
$a_{m_j} w$ to $k_c  { \phi_b  \over  c + K_c}  w$
is established by Lebesgue's dominated convergence 
theorem \cite{brezis, lions2}. The limit $\tilde c$
inherits  the same uniform upper bound on $\|\tilde c \|_{H}.$

{\it Regularity.} When $\Omega$  is a $C^2$ domain, elliptic regularity
theory implies $\tilde c^{(m)} \in H^2(\Omega)$. Sobolev embeddings 
imply $H^2(\Omega) \subset L^q(\Omega)$ for ${1/2}-{2/n}=(n-4)/2<1/q$
\cite{brezis}. If $n=2,3$ that is the case for all finite $q$ and for $q=\infty$.
$\|\tilde c^{(m)} \|_{H^2}$ and $\|\tilde c^{(m)} \|_{L^\infty}$ are
again bounded in terms of the $L^2$ norms of $\mathbf v_l$
and $\phi_b$. 

The limit $\tilde c$ inherits such $H^2$ and $L^\infty$ bounds.
Moreover, since $H^2(\Omega) \subset C(\overline \Omega)$,
$\tilde c$ must attain a minimum and a maximum value in 
$\overline \Omega$ and $\partial \Omega$.

{\it Positivity. } 
Let us now check that $c^{(m+1)}\geq 0$ by induction.
The function $c^{(1)}$ satisfies
\begin{eqnarray*}
d \int_{\Omega}  |\nabla c^{(1)-}|^2 \, d \mathbf x 
- \int_{\Omega} \mathbf v_l \cdot c^{(1)-} \nabla c^{(1)-}  \, d \mathbf x 
+ \int_{\partial \Omega^+} |c^{(1)-}|^2 \mathbf v_l \cdot
\mathbf n \, dS_\mathbf x  \\
+ \int_{\Omega} a_0 |c^{(1)-}|^2 \, d \mathbf x =0.
\end{eqnarray*}
The bilinear form $
d \int_{\Omega}  \nabla c  \nabla w \, d \mathbf x 
+ \int_{\Omega} a_0 c w \, d \mathbf x$ is coercive in $H^1(\Omega)$
when $a_0$ (that is, $\phi_b$) is not identically zero in a set of 
positive measure, which is the case. Denoting by $\lambda(d,a_0)$
this coercitivity constant, we can take $\mathbf v_l$ small enough to 
ensure $\lambda(d,a_0)/2  > k_{\mathbf v_l, \Omega} $ and get
global coercivity in $H^1$:
\begin{eqnarray*}
0  \geq (\lambda(d,a_0)- k_{\mathbf v_l, \Omega}) 
\|\tilde  c^{(1)-} \|_{H^1}^2 \geq {1 \over 2} \lambda(d,a_0)
\|\tilde  c^{(1)-} \|_{H^1}^2.
\end{eqnarray*}
This implies $c^{(1)-}=0$ and $c^{(1)}\geq 0$  provided 
$d$ is large enough and $\| \mathbf v_l  \|_{1,\infty}$ 
small enough. The same argument works for
$m \geq 1$. Notice that $a_m =  k_c  { \phi_b  \over  c^{(m)}  + K_c}
\geq k_c  { \phi_b  \over  \|c^{(m)}\|_{L^\infty}  + K_c}
\geq k_c  { \phi_b  \over  c_{\rm max} + K_c} = a_{\rm min}$, where
$ c_{\rm max} $ is a uniform upper bound of the sequence
$\|c^{(m)}\|_{L^\infty} .$ Then we can work with the bilinear
form $
d \int_{\Omega}  \nabla c  \nabla w \, d \mathbf x 
+ \int_{\Omega} a_{\rm min} c w \, d \mathbf x$ for all the
iterates, and use a uniform coeercitivity constant 
$\lambda(d,a_{\rm min})$ for all.

On one hand, this justifies the hypothesis $a_m \geq 0$ 
and $a_m  \in  L^{\infty}(\Omega)$ in the iterative scheme, so that
a solution $c \in H^1(\Omega)$ of (\ref{nutrient_st}) indeed exists. 
On the other, pointwise convergence implies $c \geq 0$.

{\it Uniqueness.}
Finally, for uniqueness assume problem (\ref{nutrient_st})
has two positive solutions $c_1$ and $c_2$ in $H^1(\Omega)$
and set $c = c_1 - c_2$. Then $u$ is a solution of 
\begin{eqnarray*}
\begin{array}{rcll}
-  d \Delta c + {\rm div}  (\mathbf v_l c)   &=&  
-  k_c K_c \phi_b { c \over (c_1 + K_c)(c_2 + K_c)}, &
  \mathbf x \in \Omega, \\
  c &=& 0, &  \mathbf x \in \partial \Omega^-, \\
  {\partial c \over \partial \mathbf n} &=& 0, &  
  \mathbf x \in \partial \Omega^+. \\
\end{array}  
\end{eqnarray*}
We write the equation in variational form and choose as test
function $w=c$.  Using the coercivity on $H$ of the part of the 
bilinear form not involving  $\phi_b$, as well as the positivity of  
$\phi_b, c_1, c_2$ and the fact that $d$
is sufficiently large and $\mathbf v_l$ sufficiently small, we
conclude that $\| c \|_H \leq 0.$ This implies
$c_1=c_2$ and the positive solution is unique.

{\it Strict positivity.} Let $c_{\rm max}=c_{\rm max}(\Omega,d,
c_0, k_c,K_c)$ an upper bound of $c$, that is,
$c \leq c_{\rm max}$ in $\Omega$.
Assuming ${\rm div}(\mathbf v_l) \leq 0$ and $0\leq \phi_b\leq 1$,
we have $- d \Delta c + \mathbf v_l \cdot \nabla c = 
-k_c {c \over c+K_c} \phi_b - {\rm div}(\mathbf v_l) c \geq 
- {k_c\over K_c}c \geq - {k_c\over K_c}c_{\rm max}$. 
By comparison principles for elliptic equations,
$c$ is bounded from below by  
$\underline{c}$ and $\underline{\underline{c}}+c_0$,
given by solutions of 
\begin{eqnarray*}
\begin{array}{rcll}
-  d \Delta \underline{c}  + \mathbf v_l \cdot \nabla \underline{c}   
+ {k_c \over K_c}  \underline{c} &=& 0,   \\
\underline{c} &=& c_0, &   \\
{\partial \underline{c} \over \partial \mathbf n} &=& 0,  
\end{array}  
\begin{array}{rclll}
- d \Delta \underline{\underline{c}} + \mathbf v_l \cdot \nabla 
\underline{\underline{c}} &=&  - {k_c\over K_c}c_{\rm max},  
&  & \mathbf x \in \Omega, \\
\underline{\underline{c}} &=& 0, & &  \mathbf x \in \partial \Omega^-, \\
{\partial \underline{\underline{c}} \over \partial \mathbf n} &=& 0, &  & 
  \mathbf x \in \partial \Omega^+,
\end{array}
\end{eqnarray*}
that is, $c \geq \underline{c} \geq \underline{\underline{c}}+c_0$.
Since $\mathbf v_l \in C(\overline \Omega)$,
$\underline{c}, \underline{\underline{c}} \in C^2(\Omega) \cap  
C(\overline \Omega)$. 
If a strict minimum value is attained at $\mathbf x \in \Omega$, 
$ d \Delta \underline{c} = {k_c \over K_c} \underline{c} >0$ at it.
The system
\begin{eqnarray*}
\begin{array}{rcll}
- d \Delta w &=&  - {k_c\over K_c}c_{\rm max},  & \mathbf x \in \Omega, \\
w & = & 0, &  \mathbf x \in \partial \Omega^-, \\
{\partial w \over \partial \mathbf n} &=& 0, &  
  \mathbf x \in \partial \Omega^+, \\
\end{array} 
\end{eqnarray*}
admits subsolutions of the form $\underline{w}(\mathbf x) = 
{k_c c_{\rm max} \over d K_c} x_n ({x_n \over 2}-M)$, $x_n$ being the 
component in the normal direction to $\partial \Omega^-$. They satisfy the 
equation and the Dirichlet condition. For $M$ large enough, depending on 
$\partial \Omega^+$, ${\partial \underline{w} \over \partial \mathbf n}
< 0$ on $\partial \Omega^+$. This requires ${k_c c_{\rm max} \over d K_c} 
(x_n -M) <0$, that is, $M$ larger than the largest vertical diameter in $\Omega$.
Under this assumption, $\underline{w}-w$ satisfies the equation with a 
zero right hand side, zero Dirichlet condition and strictly negative Neumann 
boundary condition, so that $c_0 + w \geq c_0 + \underline{w} \geq 
c_0 + \underline{w}(M) > 0 $ in $\overline{\Omega}$ provided 
$- {k_c c_{\rm max} \over 2 d K_c} M + c_0 >0$, that is, $M <  
{2d K_c \over k_c c_{\rm max}} c_0$.

Now, the function $u = w - \underline{\underline{c}}$ satisfies 
\begin{eqnarray*}
\begin{array}{rcll}
- d \Delta u + \mathbf v_l \cdot \nabla u &=&  
v_{l,n} {k_c c_{\rm max} \over d K_c} (x_n -M),  & \mathbf x \in \Omega, \\
u & = & 0, &  \mathbf x \in \partial \Omega^-, \\
{\partial u \over \partial \mathbf n} &=& 0, &  
  \mathbf x \in \partial \Omega^+. \\
\end{array} 
\end{eqnarray*}
The right hand side is negative, $u<0 $ and $c \geq \underline{\underline{c}} + c_0
\geq w + c_0 > - {k_c c_{\rm max}
\over 2 d K_c} M + c_0 >0$ in $\overline{\Omega}$.
$\square$

\section{Well-posedness results for the coupled stationary system}
\label{sec:coupled}

Consider now the full system (\ref{vfractionl_intro})-(\ref{nutrient_intro})
with boundary conditions (\ref{bc_velocity+})-(\ref{bc_concentration+})
set on a fixed domain $\Omega$.
As said before, we make the standard choice $\pi(\phi_b) = \Pi \phi_b$ in 
(\ref{velocityb_intro}). 
Numerical  simulations and asymptotic solutions \cite{entropy, seminara}
suggest the following simplification of the model:
\begin{itemize}
\item Set $\phi_l = \phi_\infty \in (0,1)$ in (\ref{pressure_intro})
and (\ref{velocityl_intro}) so that $\xi(\phi_l)
= \phi_\infty^2 {\zeta_\infty \over \mu_l} = \xi_\infty >0$
and $\eta(\phi_l) = \phi_\infty {\zeta_\infty \over \mu_l} = 
{\xi_\infty \over \phi_\infty}> 0$.
This avoids technical problems due to the low regularity of the
solutions $\phi_b$ of the stationary transport equations.
\item Set ${c \over c+K_b} = g_\infty >0$, a reference constant
term representing nutrient consumption in (\ref{nutrient_intro}).
This condition can be removed if we have  uniform positive 
lower bounds for solutions $c$ of (\ref{nutrient_intro}) for a range
of velocities $\mathbf v_l$.  Theorem 5.1 shows that this is possible 
when the height of $\Omega$ is small enough depending on $k_c$,  
$K_c$, $d$ and $c_0$, ${\rm div}(\mathbf v_l) \leq 0$, and 
$\mathbf v_{l}$  is small enough.
\end{itemize}
In Section \ref{sec:transport} we have also assumed
${\rm div}(\mathbf v_l) \leq 0$, which is associated
with $\mathbf v_l \cdot \mathbf n\leq 0$ from a physical point of view.
Notice that ${\rm div}(\mathbf v_l) \phi_l + \mathbf v_l \nabla \phi_l 
= - {k_b c \over c+K_b} \phi_b$, where $\phi_l$ is expected to be 
almost constant, that is,  $\nabla \phi_l $ is almost zero. This gives
${\rm div}(\mathbf v_l) \phi_l \sim - {k_b c \over
c+K_b} \phi_b\leq 0$, since $\phi_b \geq 0$ and $c \geq 0$.

The existence proof for the simplified system relies on an iterative scheme, 
initialized as follows. We set constant $\phi_l^{(0)} = \phi_\infty \in (0,1),$ 
$\phi_b^{(0)}= 1 - \phi_l^{(0)}$. For $\ell>0$,  the iterates are defined by 
the scheme
\begin{eqnarray}  \begin{array}{l}
\mu_b \Delta \mathbf v_b^{(\ell)} + {\mu_b \over 3}  \nabla {\rm div} 
(\mathbf v_b^{(\ell)} )  = \nabla (\Pi \phi_b^{(\ell-1)}+ p^{(\ell)}),  \\[1.5ex]
\xi_\infty \Delta p^{(\ell)}  = {\rm div}(\mathbf v_b^{(\ell)} ), 
\\[1.5ex]
\mathbf v_l^{(\ell)} = \mathbf v_b^{(\ell)} - {\xi_\infty \over
\phi_\infty}  \nabla p^{(\ell)}, 
\\[1.5ex]
{\rm div}(\mathbf v_l^{(\ell)}  \phi_l^{(\ell)} )
- k_b g_\infty  \phi_l^{(\ell)}  =  
- k_b g_\infty,  \\[1.5ex]
\phi_b^{(\ell)} = 1 - \phi_l^{(\ell)}, \\[1.5ex]
d \Delta c^{(\ell)}  - {\rm div}  (\mathbf v_l^{(\ell)}  c^{(\ell)})  
 =  k_c  { \phi_b^{(\ell)}  \over c^{(\ell)}  + K_c} c^{(\ell)},  
\end{array} \label{iterate} \end{eqnarray}
with boundary conditions (\ref{bc_velocity+})-(\ref{bc_concentration+}).
Section \ref{sec:stokes} constructs $p^{(\ell)}$ and $\mathbf v_b^{(\ell)}$.
Then, $\phi_l^{(\ell)}$ is given by section \ref{sec:transport} and $c^{(\ell)}$
by section \ref{sec:crd}.

We have the following convergence result for this scheme, leading
to a solution satisfying a number of stability bounds. \\

{\bf Theorem 6.1.} {\it Let $\Omega \subset \mathbb R^n$, 
$n=2,3,$ be a bounded domain with  $C^3$ boundary. Consider
positive constants $k_b, \mu_b,  \Pi, \xi_\infty, k_c, K_c, d, c_0$
and constant  boundary data $ \mathbf t_{\rm ext}$ and
$p_{\rm ext} - \pi_{\rm ext}.$
Provided  $\Pi, \mathbf t_{\rm ext}, p_{\rm ext} - \pi_{\rm ext}$ are 
small enough, the iterative scheme  (\ref{iterate}) converges to a 
solution 
$\mathbf v_b \in [H^2(\Omega)\cap W^{2,q}(\Omega)]^n$, $n<q<\infty$,
$p \in H^3(\Omega)\cap W^{3,q}(\Omega)$, 
$\mathbf v_l \in [H^2(\Omega)\cap W^{2,q}(\Omega)]^n$, 
$\phi_l \in H^1(\Omega) \cap W^{1,q}(\Omega)$,
$\phi_b \in  H^1(\Omega) \cap W^{1,q}(\Omega)$,
$c \in H^2(\Omega)$
of system (\ref{vfractionl_stinfty}), (\ref{velocity_st}), (\ref{pressure_st}), 
(\ref{nutrient_st}) satisfying the relations $\phi_b=1-\phi_l$ and 
$\mathbf v_l = \mathbf v_b - {\xi_\infty\over \phi_\infty} \nabla p$, as well 
as the boundary conditions (\ref{bc_velocity+})-(\ref{bc_concentration+})
on $\partial \Omega$.

Moreover, $\phi_l, \phi_b$ and $c$ are positive functions and 
the following estimates hold
\begin{eqnarray*}
0 \leq \phi_b = 1 - \phi_l \leq 1, \\
{1\over 2 }\|\mathbf v_l \|_{[W^{2,q}]^n} \leq  \left[1 + {K_{3q}(\Omega)
 \over \phi_\infty } \right] K_{2q}(\Omega,\mu_b, \xi_\infty) 
\big[  \|p_0\|_{W^{2,q}}  +  
 \| \mathbf t_{\rm ext}\|_{[L^q(\partial \Omega)]^n} \big] \\
+ {\xi_\infty \over \phi_\infty} (K_{3q}(\Omega)+1) \|p_0\|_{W^{3,q}}
= A(\Omega, \mu_b, \xi_\infty, \phi_\infty, p_0, \mathbf  t_{\rm ext}) := A, \\
 k_b g_\infty \|\nabla \phi_l \|_{L^q} \leq 4 A, \\ 
\|\mathbf v_b  \|_{[W^{2,q}]^n}  \leq 
K_{2q}(\Omega,\mu_b,\xi_\infty) \left[ {4 \Pi \over k_b g_\infty} A +
\| p_0\|_{W^{2,q}} + \|\mathbf t_{\rm ext}\|_{[L^q(\partial \Omega)]^n} \right]:= B, \\
\| p \|_{[W^{3,q}]^n} \leq {K_{3q}(\Omega) \over \xi_\infty}
B  + (K_{3q}(\Omega)+1)  \| p_0\|_{[W^{3,q}]^n},\\
{\rm min}(d,P_{\partial \Omega^-}) \|c - c_0\|_{H^1} 
\leq 4 c_0  \,  \left[2 A + {k_c \over K_c} \right],
\end{eqnarray*}
for $n< q < \infty$ and similar estimates for $q=2$.
}

{\bf Proof. }
{\it Existence of iterates and estimates. } The data $p_{\rm ext}-\pi_{\rm ext}$
satisfy all the regularity hypotheses in Theorems 4.1-4.3. Furthermore, 
$0\leq \phi_b^{(0)} \leq  1$ is constant, thus, it satisfies all the regularity 
hypotheses too.
Assume  that $\phi_b^{(\ell-1)}\in H^1\Omega)\cap W^{1,q}(\Omega)$,
$n<q<\infty$ and $0 \leq \phi_b^{(\ell-1)} \leq 1$. Then, by Theorems 4.1-4.3
we have a unique solution
$\mathbf v_b^{(\ell)} \in [W^{2,q}(\Omega) \cap H^2(\Omega)]^n$, 
$p^{(\ell)} \in W^{3,q}(\Omega) \cap H^3(\Omega)$ for the
Stokes problem. By Sobolev injections, 
$\mathbf v_b^{(\ell)}  \in [C(\overline{\Omega})]^n$,
$\mathbf v_b^{(\ell)} \in [W^{1,\infty}(\Omega)]^n$,
$p \in C^1(\overline{\Omega})$ and
$p^{(\ell)} \in W^{2,\infty}(\Omega)$.
Moreover, we have the estimates (\ref{bound_vp4})-(\ref{bound_p4}).

Using them, we define 
$\mathbf v_l^{(\ell)} \in [W^{2,q}(\Omega) \cap H^2(\Omega)]^n$,
satisfying $\mathbf v_l^{(\ell)}  \in [C(\overline{\Omega})]^n$ and
$\mathbf v_l^{(\ell)} \in [W^{1,\infty}(\Omega)]^n$.

Let us assume by now that ${\rm div}(\mathbf v_l)^{(\ell)} \leq 0$ and 
$\mathbf v_l^{(\ell)} \cdot \mathbf n \leq 0$. We also assume that 
$\|\mathbf v_l^{(\ell)}\|_{[W^{1,\infty}]^n}$ is small enough compared to 
$k_n g_\infty$. We will check these two points later in the proof.
Now, we apply Corollary 3.6 to the transport problem to 
construct a solution  $\phi_l^{(\ell)} \in H^1(\Omega) \cap W^{1,q}(\Omega)$ 
satisying $0 \leq \phi_l^{(\ell)} \leq 1.$ Next, we define 
$\phi_b^{(\ell)}=1-\phi_l^{(\ell)}  \in H^1(\Omega) \cap W^{1,q}(\Omega)$
satisfying also $0 \leq \phi_b^{(\ell)} \leq 1$. 

With the smallness assumption on $\mathbf v_l^{(\ell)}$ just made,
Theorem 5.1 provides a nonnegative solution $c^{(\ell)} \in H^2(\Omega) 
\cap C(\overline{\Omega}),$ bounded in terms of $\phi_b^{(\ell)}$ and 
$\mathbf v_l^{(\ell)}$.

{\it Uniform bounds.}
Let us denote $V^{(\ell)} = {\rm max}_{0\leq k \leq \ell} \|v_l^{(k)}\|_{[W^{2,q}]^n}$. 
Recall that $ p^{(\ell)} =  \tilde p^{(\ell)} + p_0$, $p_0 \in W^{3,q} \cap H^3$.
Theorem 4.3 guarantees estimates (\ref{bound_vp4}) and 
(\ref{bound_p4}) on $\mathbf v_b^{(\ell)}$ and  $p^{(\ell)}$.  
Using the definition of $\mathbf v_l^{(\ell)}$ we find
\begin{eqnarray*}
\| \mathbf v_l^{(\ell)} \|_{[W^{2,q}]^n} \leq \| \mathbf v_b^{(\ell)} \|_{[W^{2,q}]^n}
+ {\xi_\infty \over \phi_\infty}  \| p^{(\ell)} \|_{W^{3,q}} \leq  \\
 \left(1+{K_{3q}(\Omega)  \over \phi_\infty} \right)
 \| \mathbf v_b^{(\ell)} \|_{[W^{2,q}]^n} +  {\xi_\infty \over \phi_\infty} 
 (K_{3q}(\Omega) +1) \| p_0 \|_{W^{3,q}}.
\end{eqnarray*}
Estimate (\ref{bound_vp4}) implies then
\begin{eqnarray*}
\| \mathbf v_l^{(\ell)} \|_{[W^{2,q}]^n} \leq 
  \left[1+{K_{3q}(\Omega) \over \phi_\infty} \right]
 K_{2q}(\Omega,\mu_b,\xi_\infty)
 \Big[  \Pi \| \nabla\phi_b^{(\ell-1)} \|_{L^q}  + \| p_0\|_{W^{2,q}} \\ 
+ \| \mathbf t_{\rm ext} \|_{[L^q(\partial \Omega)]^n} \Big]
+ {\xi_\infty \over \phi_\infty}  (K_{3q}(\Omega) +1) \| p_0} \|_{W^{3,q}.
\end{eqnarray*}
Notice that $\nabla \phi_b^{(\ell-1)} = - \nabla \phi_l^{(\ell-1)}$.
Now, Corollary 3.6. under the smallness assumption
ensures that
\begin{eqnarray*}
{1\over 2} k_b g_\infty \| \nabla \phi_l^{(\ell-1)} \|_{L^2} \leq 
\| \nabla {\rm div}(\mathbf v_l^{(\ell-1)}) \|_{[L^q]^n} \leq
\| \mathbf v_l^{(\ell-1)} \|_{[W^{2,q}]^n}. 
\end{eqnarray*}
Combining these inequalities, we find
\begin{eqnarray*}
V^{(\ell)}  \leq     \left[1+{K_{3q}(\Omega) \over \phi_\infty} \right]
 K_{2q}(\Omega,\mu_b,\xi_\infty)
 \Big[  {2 \Pi \over k_b g_\infty} V^{(\ell)} + \| p_0\|_{W^{2,q}} \\ 
 + \| \mathbf t_{\rm ext} \|_{[L^q(\partial \Omega)]^n} \Big]
 + {\xi_\infty \over \phi_\infty} (K_{3q}(\Omega) +1) \| p_0} \|_{W^{3,q}.
\end{eqnarray*}
If   $\left[1+{K_{3q}(\Omega) \over \phi_\infty} \right]
K_{2q}(\Omega,\mu_b,\xi_\infty) {2 \Pi \over k_b g_\infty} < 1/2$, then
$V^{(\ell)} $ is uniformly bounded by a right hand side which
does not depend on $\ell$. Since 
$0 \leq \phi_l^{(\ell)} = 1 - \phi_b^{(\ell)} \leq 1$, the uniform
bound on $\|v_l^{(\ell)}\|_{[W^{2,q}]^n}$ extends to 
$\|\phi_l^{(\ell)} \|_{W^{1,q}}$ and $\|\phi_b^{(\ell)} \|_{W^{1,q}}$.
It also extends  $\mathbf v_b^{(\ell)}$ and $p^{(\ell-1)}$ in view of
\begin{eqnarray*}
\|\mathbf v_b^{(\ell)} \|_{[W^{2,q}]^n}  \leq 
K_{2q}(\Omega,\mu_b,\xi_\infty) \left[ \Pi \| \nabla \phi_b^{\ell-1}\|_{L^q} +
\| p_0\|_{W^{2,q}} + \|\mathbf t_{\rm ext}\|_{[L^q(\partial \Omega)]^n} \right], \\
\| p^{(\ell)} \|_{[W^{3,q}]^n} \leq {K_{3q}(\Omega) \over \xi_\infty}
\| \mathbf v_b^{(\ell)} \|_{[W^{2,q}]^n}  +
 (K_{3q}(\Omega) +1) \| p_0\|_{[W^{3,q}]^n}.
\end{eqnarray*}
and to $\|c^{(\ell)} \|_{H^2}$ by Theorem 5.1

{\it Smallness and sign constraints.} The previous two steps require 
smallness of $\|\mathbf v_l^{(\ell)}|_{[W^{1,\infty}]}$ and 
${\rm div}(\mathbf v_l^{(\ell)}) \leq 0$, $\mathbf v_l^{(\ell)} \cdot 
\mathbf n \leq 0$ to be able to apply Corollary 3.6 and Theorem 5.1.  
Let us proceed by induction to guarantee these conditions.

Initially, $\phi_b^{(0)}$ is constant, thus $\nabla \phi_b^{(0)}=0$. We 
can apply Theorems 4.1-4.3 to construct $\mathbf v_b^{(1)}$ and 
$p^{(1)}$ in such a way that $\|\mathbf v_b^{(1)} \|_{[W^{2,q}]^n}$,
$\|p^{(1)} \|_{[W^{3,q}]^n}$ and  $\|\mathbf v_l^{(1)} \|_{[W^{2,q}]^n}$ 
are bounded in terms of the problem parameters $\|p_0\|_{[W^{3,q}]^n}$ 
and $\|\mathbf t_{\rm ext}\|_{[L^q(\partial \Omega)]^n}$. By Sobolev
injections with $n < q < \infty$, $\|\mathbf v_b^{(1)} \|_{[W^{1,\infty}]^n}$ 
satisfies  a similar estimate, and can be made as small as needed by 
making $\mathbf t_{\rm ext}$ and $p_{\rm ext}-\pi_{\rm ext}$ small.
By Corollary 3.6, $\|\nabla \phi_l^{(1)}\|_{L^q}$   is then bounded by 
$\|\mathbf v_l^{(1)} \|_{[W^{2,q}]^n}$  and is equally small.

Moreover, ${\rm div}(\mathbf v_l^{(1)}) \phi_l^{(1)} + \mathbf v_l^{(1)} 
\nabla \phi_l^{(1)} = - k_b g_\infty \phi_b^{(1)} \leq 0$. Both
$\mathbf v_l^{(1)} $ and $\nabla \phi_l^{(1)} $ are small compared
to $\phi_l^{(1)}$ and $- k_b g_\infty \phi_b^{(1)}$ which are 
almost constant. Thus, ${\rm div}(\mathbf v_l^{(1)}) \leq 0$.
Now, $\int_{A} {\rm div}(\mathbf v_l^{(1)}) \, d {\mathbf x}
= \int_{\partial A} \mathbf v_l^{(1)} \cdot \mathbf n \, d {S_{\mathbf x}}
\leq 0$ for any $A \subset \Omega$ implies
$\mathbf v_l^{(1)} \cdot \mathbf n \leq 0$ on $\partial \Omega$.

By induction, assuming that $\|\mathbf v_b^{(\ell-1)} \|_{[W^{1,\infty}]^n}$
satisfies the smallness requirement, we can repeat the same procedure
to prove that $\|\mathbf v_b^{(\ell)} \|_{[W^{1,\infty}]^n}$ satisfies it
too and that it also satisfies the sign conditions. The only difference is that
now we need to estimate $\| \nabla {\rm div}(\mathbf v_l^{(\ell-1)}) \|_{[L^q]^n}$
and require that $\Pi$ is small enough too.

{\it Convergence to a solution.} The uniform bounds we
have just established allow us to prove convergence of
the iterative scheme to a solution. 
We can extract subsequences converging to limits $\phi_l=1-\phi_b$ 
weakly in $H^1$ and strongly in $L^2$, limits $v_{l,j}$,
 $v_{b,j}$ $j=1,\ldots,n$, and $p$  weakly in $W^{2,q}\cap H^2$, 
$n<q<\infty$, strongly in $W^{1,q}\cap L^2$ and strongly in 
$C(\overline{\Omega})$. We have
\begin{eqnarray*}
\int_{\Omega}  \phi_l^{(\ell)} \mathbf v_l^{(\ell)}
\cdot \nabla w \, d \mathbf x 
- \int_{\partial \Omega}  \phi_l^{(\ell)} \mathbf v_l^{(\ell)} 
\cdot \mathbf n dS_{\mathbf x} \\
+ \int_{\Omega} k_b g_\infty \phi_l^{(\ell)} w \, d \mathbf x
= \int_{\Omega} k_b g_\infty w \, d \mathbf x.
\end{eqnarray*}
Strong convergence of $\mathbf v_l^{(\ell)}
\cdot \mathbf n$ in $C(\partial \Omega)$
and $ \mathbf v_l^{(\ell)}   \cdot \nabla w$
in $L^2(\Omega)$ and weak convergence
of $\phi_l^{(\ell)}$ allow us to pass to the
limit as $\ell \rightarrow \infty$ and prove that
$\phi$ is a solution of the stationary transport 
problem for $\mathbf v_l$, given $w \in H^1(\Omega)$.

Next, given $(\mathbf w, q) \in [H^1_{0,-}(\Omega)]^n \times 
H^1_0(\Omega)$ we have
\begin{eqnarray*} 
 \int_{\Omega}
[ \mu_b \nabla \mathbf v_b^{(\ell)}  \nabla \mathbf w  + {\mu_b \over 3} 
{\rm div}(\mathbf v_b^{(\ell)}) {\rm div} (\mathbf w ) ] \, d \mathbf x 
- \int_{\Omega}   \tilde p^{(\ell)}  \, {\rm div}(\mathbf w) \, d \mathbf x 
+ \int_{\Omega} \xi_\infty \nabla \tilde p^{(\ell)} \nabla q \, d \mathbf x \\
+ \int_{\Omega} {\rm div}(\mathbf v_b^{(\ell)} ) q \, d \mathbf x
= \int_{\Omega}  \Pi \phi_b^{(\ell-1)} \,{\rm div} (\mathbf w) \, d \mathbf x
 + \int_{\partial \Omega}   \mathbf t_{\rm ext} \mathbf w \, d S_{\mathbf x }
- \int_{\Omega} \xi_\infty\nabla p_0 \nabla q \, d \mathbf x.
 \end{eqnarray*}
Weak $L^2$ convergence is enough to pass to the limit
in all terms and  find
\begin{eqnarray*} 
 \int_{\Omega}
[ \mu_b \nabla \mathbf v_b  \nabla \mathbf w  + {\mu_b \over 3} 
{\rm div}(\mathbf v_b) {\rm div} (\mathbf w ) ] \, d \mathbf x 
- \int_{\Omega}   \tilde p  \, {\rm div}(\mathbf w) \, d \mathbf x 
+ \int_{\Omega} \xi_\infty \nabla \tilde p \nabla q \, d \mathbf x \\
+ \int_{\Omega} {\rm div}(\mathbf v_b) q \, d \mathbf x
= \int_{\Omega}  \Pi \phi_b \,{\rm div} (\mathbf w) \, d \mathbf x
 + \int_{\partial \Omega}   \mathbf t_{\rm ext} \mathbf w \, d S_{\mathbf x }
 - \int_{\Omega} \xi_\infty\nabla p_0 \nabla q \, d \mathbf x.
 \end{eqnarray*}
We set $p = \tilde p + p_0$. Notice that the weak limits
satisfy $(\mathbf v_b, \tilde p) \in [H^1_{0,-}(\Omega)]^n
\times  H^1_0(\Omega)$.

Finally, given $w \in H^1_{0,-}(\Omega)$ we have
 \begin{eqnarray*}
d \int_{\Omega}  \nabla \tilde c^{(\ell)} 
 \nabla w \, d \mathbf x 
- \int_{\Omega} \mathbf v_l^{(\ell)} \tilde c^{(\ell)} \nabla w
 \, d \mathbf x 
+ \int_{\partial \Omega} \tilde c^{(\ell)} w \mathbf v_l^{(\ell)} \cdot
\mathbf n \, d {S_\mathbf x}  \\
 + \int_{\Omega} {k_c \phi_b^{(\ell)} \over c^{(\ell)}  + K_c} 
 \tilde c^{(\ell)} w \, d \mathbf x    
 = c_0 \int_{\Omega} \mathbf v_l^{(\ell)}  \nabla w \, d \mathbf x   
 - c_0\int_{\Omega} {k_c \phi_b^{(\ell)}\over c^{(\ell)}  + K_c}  w \, 
d \mathbf x,
\end{eqnarray*}
for  $\tilde c^{(\ell)}$ uniformly bounded in $H^2\cap H^1_{0,-}$.
We can extract a subsequence converging to a limit
$\tilde c \in H^2\cap H^1_{0,-}$, weakly in $H^2$ and strongly, 
at least in $L^2$ and $L^4$, as well as pointwise in $\Omega$. 
The traces of $c^{(\ell)}$ in $\partial \Omega$
converge weakly in $L^2(\partial \Omega)$, while 
$\mathbf v_l^{(\ell)} w $ converges strongly in 
$L^2(\partial \Omega)$, as argued above. This allows
us to pass to the limit in all the integrals except those
involving $\phi_b^{(\ell)}.$ For those, we remark that
the integrands converge pointwise and they are uniformly
bounded by $L^\infty$ functions. Convergence follows
by Lebesgue's Theorem.
Taking limits in the identity we find
 \begin{eqnarray*}
d \int_{\Omega}  \nabla \tilde c  \nabla w \, d \mathbf x 
- \int_{\Omega} \mathbf v_l \tilde c  \nabla w
 \, d \mathbf x 
+ \int_{\partial \Omega} \tilde c w \mathbf v_l  \cdot
\mathbf n \, d {S_\mathbf x}  \\
 + \int_{\Omega} {k_c \phi_b  \over c   + K_c} 
  \tilde c w \, d \mathbf x    
 = - c_0 \int_{\Omega} \mathbf v_l \nabla w \, d \mathbf x   
 - c_0\int_{\Omega} {k_c \phi_b \over c + K_c}  w \, 
d \mathbf x,
\end{eqnarray*}
and set $c = \tilde c + c_0$. The functions $\mathbf v_l$,
$\mathbf v_b$, $p$, $\phi_l\geq 0$, $\phi_b=1-\phi_l \geq 0$, 
and $c \geq 0$ provide the solution we seek. 

Passing to the limit on the uniform bounds established for
the convergent sequences, the same stability
bounds hold for the solution.
$\square$ \\

 
In a similar way, we can handle the model with the original
stationary transport problem. \\

{\bf Theorem 6.2.} {\it We keep the hypotheses made in Theorem
6.2 and consider the system (\ref{vfractionl_st}), (\ref{velocity_st}), 
(\ref{pressure_st}),  (\ref{nutrient_st}) with the additional  relations 
$\phi_b=1-\phi_l$ and $\mathbf v_l = \mathbf v_b - 
{\xi_\infty\over \phi_\infty} \nabla p$, and
the boundary conditions (\ref{bc_velocity+})-(\ref{bc_concentration+}).
Assume that $M= {\rm max} \{ x_n \, | \, (x_1,\ldots,x_{n-1},x_n) \in 
\Omega\} < {2d K_c \over k_c c_{\rm max}} c_0$, 
we have a solution
$\mathbf v_b \in [H^2(\Omega) \cap W^{2,q}(\Omega)]^n$,
$p \in H^3(\Omega)\cap W^{3,q}(\Omega)$, 
$\mathbf v_l \in [H^2(\Omega) \cap W^{2,q}(\Omega)]^n$, 
$\phi_l \in H^1(\Omega) \cap W^{1,q}(\Omega)$,  
$\phi_b \in H^1(\Omega) \cap W^{1,q}(\Omega)$, 
$c\in H^2(\Omega)$
of system (\ref{vfractionl_st}), (\ref{velocity_st}), (\ref{pressure_st}), 
(\ref{nutrient_st}) 
satisfying the relations $\phi_b=1-\phi_l$ and $\mathbf v_l = \mathbf v_b - 
{\xi_\infty\over \phi_\infty} \nabla p$, as well as 
the boundary conditions (\ref{bc_velocity+})-(\ref{bc_concentration+}).
This solution has the same regularity as the one obtained in Theorem
6.1 and satisfies the same estimates, replacing $k_b g_\infty$ with
$ {c_{\rm min} k_b \over c_{\rm max} + K_b}$.
}

{\bf Proof.}
For $\ell>0$,  we consider the iterative scheme
\begin{eqnarray}  \begin{array}{l}
\mu_b \Delta \mathbf v_b^{(\ell)} + {\mu_b \over 3}  \nabla {\rm div} 
(\mathbf v_b^{(\ell)} ) = \nabla (\Pi \phi_b^{(\ell-1)}+ p^{(\ell)}),  
\\[1.5ex]
\xi_\infty \Delta p^{(\ell)} = {\rm div}(\mathbf v_b^{(\ell)}), 
\\[1.5ex]
\mathbf v_l^{(\ell)} = \mathbf v_b^{(\ell)} - {\xi_\infty \over
\phi_\infty}  \nabla p^{(\ell)}, 
\\[1.5ex]
{\rm div}(\mathbf v_l^{(\ell)}  \phi_l^{(\ell)} )
- k_b {c^{\ell-1} \over c^{\ell-1} + K_b} \phi_l^{(\ell)}  =  
- k_b {c^{\ell-1} \over c^{\ell-1} + K_b}  \\ [1.5ex]
\phi_b^{(\ell)} = 1 - \phi_l^{(\ell)}, 
\\[1.5ex]
d \Delta c^{(\ell)}  - {\rm div}  (\mathbf v_l^{(\ell)}  c^{(\ell)})  
 =  k_c  { \phi_b^{(\ell)}  \over c^{(\ell)}  + K_c} c^{(\ell)},  
\end{array} \label{iterate2} \end{eqnarray}
with boundary conditions (\ref{bc_velocity+})-(\ref{bc_concentration+}),
starting from constant $\phi_l^{(0)} = \phi_\infty \in (0,1),$ 
$\phi_b^{(0)}= 1 - \phi_l^{(0)}$, $c^{(0)}=c_0$. 

The only difference with respect to the previous proof is that the constant
$g_\infty$ in the stationary transport equation is replaced by the function
$g(c) = {c \over c + K_b}$. Existence of a solution with enough
regularity properties is guaranteed  by Corollary 3.7. At each step
$c^{(\ell-1)} \in H^2(\Omega)$ is strictly positive and uniformly bounded from 
below by a value $c_{\rm min}$ for all $\ell$.  Thus, $0\leq g(c^{(\ell)})
\leq 1$.  Moreover, for  $\|\mathbf v_l^{(\ell)} \|_{[W^{1,q}]^n}$ smaller
than a constant independent of $\ell$ and keeping the notation of the
proof of Theorem 5.1, we have
\begin{eqnarray*}
\|c^{(\ell)} - c_0\|_{H^2} \leq K(\Omega,d,c_0) \,  
\left[ \| \mathbf v_l^{(\ell)}\|_{[L^2]^n} 
+ {k_c \over K_c} \|\phi_b^{(\ell)}\|_{L^2} \right],
\end{eqnarray*}
which provides uniform bounds on $\|c^{\ell}\|_{L^\infty}$ inherited by the
limit $c$.  To ensure $g(c^{(\ell)})\in W^{1,q}(\Omega)$ we need 
$c^{(\ell)} \in W^{1,q}(\Omega)$, which follows by elliptic regularity
for the Laplacian with a right hand side in $W^{-1,q}(\Omega)$. This
also provides a uniform $W^{1,q}(\Omega)$ bound on $g(c^{(\ell)})$.
$\square$

\section{Conclusions and perspectives}
\label{sec:perspectives}

We have studied a quasi-stationary system governing biofilm
spread on surfaces in terms of a two phase flow mixture.
Similar models arise in the study of cells and tissues, where 
one phase is a liquid solution, whereas the other one is assorted 
biomass. 
The system combines stationary transport equations for the volume 
fractions  of liquids and biomass, compressible Stokes type systems
and Darcy relations for velocities and pressure together with
convection-reaction-diffusion systems for nutrients.
We are able to construct solutions satisfying a number of stability 
bounds under sign assumptions on the divergence and normal
components of velocity fields, motivated by numerical and asymptotic
solutions  for thin flat films.

This study lays the basis for the analysis of the lubrication type
equations describing the motion of the biofilm boundary as
it grows. For a two dimensional slice, we would have 
nonlocal equations of the form
\begin{eqnarray*}
h_t + v_1(x,h(x)) h_x + \int_0^{h(x)} v_{1,x} (x,s) ds = v_3(x,0),
\end{eqnarray*}
or
\begin{eqnarray*}
h_t + {\partial \over \partial x} \int_0^{h(x)} v_{1} (x,s) ds = v_3(x,0),
\end{eqnarray*}
where $v_j = v_{b,j} - {\xi_\infty \over \phi_\infty }  {\partial p
\over \partial x_j}$, $j=1,2$, are defined by solutions of   quasi-stationary 
systems for varying domains with upper boundary defined by $h(x,t)$.
In higher dimensions, we have equations of the form (\ref{height2})
and (\ref{height3}).
Nonlocality and the coupling with the quasi-stationary system
render this problem much more complex than classical lubrication 
models \cite{bertozzi}. 

\vskip 1cm

{\bf Acknowledgements.} 
This research has been partially supported by the FEDER /Ministerio
de Ciencia, Innovaci\'on y Universidades - Agencia Estatal de Investigaci\'on
grant PID2020-112796RB-C21. \\

\end{document}